\documentstyle[12pt,twoside]{amsart}

\renewcommand{\c}[0]{{\mathbb C}}  

\renewcommand{\o}[0]{{\cal O}} 
\newcommand{\z}[0]{{\mathbb Z}}
\newcommand{\n}[0]{{\mathbb N}}
\renewcommand{\r}[0]{{\mathbb R}} 

\renewcommand{\a}[0]{{\mathbb A}}

\newcommand{\p}[0]{{\mathbb P}}

\newcommand{\q}[0]{{\mathbb Q}}
\newcommand{\map}[0]{\dasharrow}
\newcommand{\qtq}[1]{\quad\mbox{#1}\quad}

\newcommand{\gal}[0]{\operatorname{Gal}}

\newcommand{\mult}[0]{\operatorname{mult}}

\newcommand{\supp}[0]{\operatorname{Supp}}    
\newcommand{\red}[0]{\operatorname{red}}

\newcommand{\cent}[0]{\operatorname{center}}

\newcommand{\inter}[0]{\operatorname{Int}}    
\newcommand{\sing}[0]{\operatorname{Sing}}    
\newcommand{\ex}[0]{\operatorname{Ex}}

\newsymbol\subsetneq 2328
\newsymbol\onto 1310
\def\into{\DOTSB\lhook\joinrel\rightarrow}
\newsymbol\twoheadrightarrow 1310

\newtheorem{thm}{Theorem}[section]
\newtheorem{question}[thm]{Question}
\newtheorem{lem}[thm]{Lemma}
\newtheorem{cor}[thm]{Corollary}

\newtheorem{prop}[thm]{Proposition}

\newtheorem{conj}[thm]{Conjecture}
\newtheorem{complement}[thm]{Complement}

\theoremstyle{definition}
\newtheorem{defn}[thm]{Definition}
\newtheorem{condition}[thm]{Condition}
\newtheorem{say}[thm]{}
\newtheorem{exmp}[thm]{Example}

\newtheorem{rem}[thm]{Remark}          
\newtheorem{ack}{Acknowledgments}        
\newtheorem{notation}[thm]{Notation}

\theoremstyle{remark}

\begin{document}
\bibliographystyle{amsplain}

\title{Real Algebraic Threefolds II.\\ Minimal Model Program}
\author{J\'anos Koll\'ar}

\maketitle
\tableofcontents

\section{Introduction}

In real algebraic   geometry, one of the main directions of
investigation  is the   topological study of the set of  real solutions
of algebraic equations. The first general result was proved in
\cite{Nash52},  and later developed by many others (see
\cite{AK92} for some recent directions). One of these theorems says
that every compact differentiable manifold can be realized as the set
of real points of an algebraic variety.
\cite{Nash52} posed the problem of obtaining similar results using a
restricted class of varieties, for instance rational varieties. For
real algebraic surfaces this question was settled in 
\cite{Comessatti14}.

The aim of this series of papers is to utilize the  theory of  minimal
models  to investigate this question for real algebraic threefolds.
This approach is very similar in spirit to the one employed  by
\cite{Comessatti14}. (See \cite{Silhol89, ras}  for   introductions to
real algebraic surfaces from the point of view of the minimal model
program.)

For algebraic threefolds over $\c$, the minimal model program (MMP for
short) provides a very powerful tool.  The method of the program is
the following. (See \cite{koll87, CKM88} or \cite{KM98} for
introductions.)

Starting with a smooth projective 3-fold $X$, we perform a series of
``elementary" birational transformations
$$
 X=X_0\map X_1\map \cdots \map X_n=:X^*
$$
 until we reach a variety $X^*$ whose global structure is ``simple". 
(Neither the intermediate steps $X_i$ nor the final 
$X^*$ are uniquely determined by $X$.) In essence the minimal model
program  allows us to investigate many questions in two steps: first
study the effect of the ``elementary" transformations and then
consider the ``simple" global situation. 

In practice both of these steps are frequently rather difficult. For
instance, we still do not have a complete list of all possible
``elementary" steps, despite repeated attempts to obtain it.

 A somewhat unpleasant feature of the theory is that the varieties
$X_i$ are not smooth, but have so called terminal singularities.  This
means that   $X_i(\r)$ is not necessarily a manifold. In developing
the theory of minimal models for real algebraic threefolds, we again
have to understand the occurring terminal singularities. This was done
in the first paper of this series \cite{rat1}.

If $X$ is defined over a field $K$, then there is a variant of the MMP
where  the intermediate varieties $X_i$ are also defined over
$K$. I refer to this as the MMP over $K$. This suggests the following
two step approach to understand the topology of
$X(\r)$:
\begin{enumerate}
\item  Study the topological effect of the ``elementary"
transformations.
\item Investigate  the topology of $X^*(\r)$.
\end{enumerate}

\noindent The aim of this paper is to complete the first of these two
steps.

I am unable to say much about this question in general. There are
serious problems coming from  algebraic geometry and also   from
3-manifold topology. Some of these are discussed in section 4. My aim
is therefore   more limited: find   reasonable conditions which ensure
that the steps of the MMP can be described topologically.

The simplest case   to study is contractions
$f:X\to Y$ where $X$ is smooth. Over $\c$ the complete list of such
contractions is known \cite{Mori82}, and it is not hard to obtain a
complete list over $\r$.  From this list one can see that in all such 
examples where $X(\r)\to Y(\r)$ is  complicated,
$X(\r)$ contains  a   special surface of nonnegative Euler
characteristic.  This turns out to be a general pattern, though the
proof presented here relies on  a laborious case analysis. The
precise  technical theorem is stated in (\ref{int.nonorient.thm}).

None of the complicated examples occur if $X(\r)$ is orientable, and
this yields the following:

\begin{thm}\label{int.orient.thm}
 Let $X$ be a smooth, projective, real algebraic $3$-fold and $X^*$
the result of the MMP over $\r$. Assume that  $X(\r)$ is orientable.

Then the topological normalization  $\overline{X^*(\r)}$ of $X^*(\r)$
is a  PL-manifold, and
 $X(\r)$  can be obtained from
$\overline{X^*(\r)}$ by repeated application of the following
operations: 
\begin{enumerate}
\setcounter{enumi}{-1}
\item throwing away all isolated points of $\overline{X^*(\r)}$,
\item taking connected sums of connected components,
\item taking connected sum with $S^1\times S^2$,
\item taking connected sum with $\r\p^3$.
\end{enumerate}
\end{thm}

\begin{rem} $X^*$ uniquely determines (\ref{int.orient.thm}.0) and
also (\ref{int.orient.thm}.1). The latter can be seen by analyzing
real analytic morphisms $h:[0,1]\to X^*(\r)$ where  the endpoints map
to different connected components of 
$\overline{X^*(\r)}$. In practice this may be quite hard, and it could
be easier to work through the MMP backwards.

$X^*$ contains some information about the steps
(\ref{int.orient.thm}.2--3), but these are by no means unique. Even if
$X^*$ is smooth, both of these steps are possible, as shown by the
next example.
\end{rem}

\begin{exmp}\label{connsum.exmp}
 It is well known how to create  connected sum with $\r\p^3$
algebraically.  Let $X$ be a smooth 3-fold over $\r$ and $0\in X(\r)$ a
real point. Set $Y=B_0X$. Then
$Y(\r)\sim X(\r)\ \#\  \r\p^3$.  (The   connected sum of two
nonoriented manifolds is, in general,  not unique. It is, however,
unique if one of the summands has an automorphism with an isolated
fixed point which reverses local orientation there.)

Connected sum with $S^1\times S^2$ is somewhat harder. Let $X$ be a
smooth 3-fold over $\r$ and 
$D\subset X$  a real curve which has a unique real point
$\{0\}= D(\r)$. Assume furthermore that  near $0$ the curve is given by
equations
$(z=x^2+y^2=0)$. Set
$Y_1=B_DX$. $Y_1$ has a unique singular point $P$; set $Y=B_PY_1$. It
is not hard to see that $Y$ is smooth and 
$Y(\r)\sim X(\r)\ \#\  (S^1\times S^2)$. 
\end{exmp}

\begin{rem}\label{1.basic.top.facts}
 As (\ref{int.orient.thm}) already shows, we have to move between
topological, PL and differentiable manifolds. In dimension 3 every
compact topological 3--manifold carries a unique PL--manifold
structure (cf.\ \cite[Sec.\ 36]{Moise77}) and also a unique
differentiable structure (cf.\ \cite[p.3]{Hempel76}). I mostly use the
PL--structure since  most algebraic constructions are natural in the
PL--category. For instance, $\r^1\to \r^2$ given by $t\mapsto
(t^2,t^3)$ is a PL--embedding but not a differentiable embedding in the
natural differentiable structures.

In dimension 3 the PL--structure behaves very much like a
differentiable  structure. For instance, let $M^3$ be a PL 3--manifold,
$N$ a compact PL--manifold of  dimension 1 or 2 and $g:N\into M$ a
PL--embedding. Then a suitable open  neighborhood of $g(N)$ is
PL--homeomorphic to a real vector bundle over $N$  (cf.\ \cite[Secs.\
24 and 26]{Moise77}). (Note that a similar result fails  for
topological 3--manifolds (cf.\ \cite[Sec.\ 18]{Moise77}), and it also
fails for PL 4--manifolds: take any nontrivial knot in $S^3$ and
suspend it   in $S^4$.)
\end{rem}

\begin{say}[Surfaces in 3--manifolds]\label{surf.in.3-man}
 Let $M$ be a PL 3--manifold without boundary,
$N$ a compact PL 2--manifold without boundary and $g:N\into M$ a
PL--embedding. As we noted above, a neighborhood of $N$ is an
$\r$-bundle over $N$. $\r$-bundles over $N$ are classified by group
homomorphisms $\rho:\pi_1(N)\to \{\pm 1\}$. If $\rho $ is trivial then
$N$ is 2--sided in $M$, otherwise it is 1--sided. We also allow self
homeomorphisms of $N$, thus we get the following possibilities when
$N$ has nonnegative Euler characteristic:
\begin{description}
\item[$S^2$] Always 2--sided, many such surfaces in every $M^3$.

\item[$\r\p^2$] $M^3$ is not orientable in the  2--sided case. Such
manifolds are called $\p^2$-reducible (cf.\
\cite[p.88]{Hempel76}).  In the 1--sided case the boundary of a regular
neighborhood is $S^2$, thus $M\sim M'\ \#\ \r\p^3$ for some 3--manifold
$M'$.  Most 3--manifolds do not contain any $\r\p^2$.

\item[{\rm Torus}] The 2--sided case occurs in any 3--manifold as the
boundary of a regular neighborhood of any  $S^1$ along which $M$ is
orientable. There is a unique 1--sided case. For these $M$ is not
orientable. Most nonorientable 3--manifolds do not contain 1--sided
tori, see section 12.

\item[{\rm Klein bottle}] $M$ is nonorientable in the 2--sided case.
 The boundary of a regular neighborhood of any  $S^1$ along which $M$
is nonorientable is such. There are two different 1--sided cases,
depending on whether $M$ is orientable near $N$ or not. These are
again rare, see section 12.
\end{description}

\noindent This shows that there are many 3--manifolds which do not
contain
$\r\p^2$, 1--sided tori or  Klein bottles. These correspond to 6
different cases on the above list. It turns out that we need to exclude
only 3 of these for our main theorem.
\end{say}

\begin{condition}\label{int.no.cond}
 Let $M$ be a  PL 3-manifold without boundary. Consider the following
properties:
\begin{enumerate}
\item $M$ does not contain a 2-sided $\r\p^2$,
\item $M$ does not contain a 1-sided torus,
\item $M$ does not contain a 1-sided Klein bottle with nonorientable
neighborhood.
\end{enumerate} Failure  of any of these properties implies that $M$
is not orientable, but  there are many nonorientable 3-manifolds which
do    satisfy all 3 of the above conditions. For instance, this holds
if
$M$ is hyperbolic (\ref{hyp.doesnotcont.thm}). 
\end{condition}

\begin{thm}\label{int.nonorient.thm}
 Let $X$ be a smooth, projective, real algebraic $3$-fold and $X^*$
the result of the MMP over $\r$. Assume that  $X(\r)$   satisfies  
the 3 conditions (\ref{int.no.cond}.1--3).

Then the conclusions of (\ref{int.orient.thm}) hold.
\end{thm}

\begin{rem} It would seem that we also need to allow connected sum
with 
$S^1\tilde{\times} S^2$ (cf.\ (\ref{5.notation})), corresponding to
attaching a nonorientable 1--handle. This, however, would give a
1--sided torus which we excluded.

All 3 conditions (\ref{int.no.cond}.1--3) are necessary for the
theorem to hold. My feeling is that  essentially nothing can be said
without 
 (\ref{int.no.cond}.1) or (\ref{int.no.cond}.3).  
(\ref{int.no.cond}.2) has a twofold role in the proof. First, it
ensures that $X$ is not obtained as a blow up of a smooth 3-fold $Y$
along a curve. This in itself would not be a problem, but it may
happen that $Y(\r)$ contains a 2-sided
$\r\p^2$ but $X(\r)$ does not.  It seems to me that this leads to
rather complicated topological questions.  Still, a suitable
reformulation of the theorem may get around this problem. Second,
(\ref{int.no.cond}.2) is also used to exclude a few singularities on
the $X_i$. These cases  are of index 1 and they can be described very
explicitly. It should be possible to work with them.
\end{rem}

The technical heart of the proof is a listing of the possible
singularities that occur in the course of the MMP and a fairly detailed
description of the steps of the MMP. The final result is relatively
easy to state but the proof is  a case-by-case examination.

\begin{thm}\label{int.mmp.sings}
 Let $X$ be a smooth, projective, real algebraic $3$-fold and assume
that  $X(\r)$   satisfies   the 3 conditions (\ref{int.no.cond}.1--3).

Let $X_i$ be any of the intermediate steps of the MMP over $\r$
starting with $X$ and $0\in X_i(\r)$ a real point. Then a neighborhood
of $0\in X_i$ is real analytically equivalent to one of the following
standard forms:
\begin{enumerate}

\item ($cA_0$) Smooth point.

\item ($cA_{>0}^+$) 
 $(x^2+y^2+g_{\geq 2}(z,t)=0)$, where  $g$ is not everywhere negative
in a punctured neighborhood of $0$.

\item ($cE_6$)  
$(x^2+y^3+(z^2+t^2)^2+ yg_{\geq 4}(z,t)+g_{\geq 6}(z,t)=0)$.
\end{enumerate}
\end{thm}

\begin{rem}   The symbol $g_{\geq m}$ denotes a power series of
multiplicity at least $m$. 

The name of the cases  is explained in \cite{rat1}. 

The above points of type $cE_6$ 
 form a codimension 7 family in the space of all $cE_6$ singularities.
They all occur, even if $X(\r)$ is orientable. Points of type
$cA_{>0}^+$ occur for many choices of $g$. Section 10 gives an
algorithm to decide which cases of $g$ do occur, but I was unable to
write the condition in closed form. For the applications this does not
seem to matter.
\end{rem}

Using \cite[4.3, 4.4, 4.9]{rat1}, this immediately implies:

\begin{cor}\label{int.mmp.sings.top} Notation and assumptions as in
(\ref{int.mmp.sings}). Then
$\overline{X_i(\r)}\setminus\{\mbox{isolated points}\}$ is a compact
PL 3-manifold without boundary.\qed
\end{cor}

The next step is to understand the  ``elementary"  steps of the MMP
over $\r$. (\ref{int.nonorient.thm}) turns out to be a  consequence of
(\ref{int.mmp.steps}). (See (\ref{8.wbup})  for the definition of
weighted blow-ups.)

\begin{thm}\label{int.mmp.steps}
 Let $X$ be a smooth, projective, real algebraic $3$-fold such that 
$X(\r)$   satisfies the   conditions (\ref{int.no.cond}.1--3).

Let $f_i:X_i\map X_{i+1}$ be any of the intermediate steps of the MMP
over
$\r$ starting with $X$.   Then the induced map $f_i:X_i(\r)\to
X_{i+1}(\r)$ is everywhere defined and the following is a complete
list of possibilities for $f_i$:
\begin{enumerate}

\item ($\r$-trivial) $f_i$ is an isomorphism in a (Zariski)
neighborhood of the set of real points.

\item ($\r$-small) $f_i:X_i(\r)\to X_{i+1}(\r)$ collapses a 1-complex
to   points and there are small perturbations $\tilde f_i$ of
$f_i$ such  $\tilde f_i: \overline{X_i(\r)}\to \overline{X_{i+1}(\r)}$
is a PL-homeomorphism.

\item (smooth point blow up)  $f_i$ is the inverse of the blow up of a
smooth point $P\in X_{i+1}(\r)$.

\item (singular point blow up) $f_i$ is the inverse of a (weighted) 
blow up of a singular point $P\in X_{i+1}(\r)$. There are two cases:

\begin{enumerate}

\item ($cA_{>0}^+$, $\mult_0g$ even)  Up to real analytic equivalence
near $P$,
 $X_{i+1}\cong (x^2+y^2+g_{\geq 2m}(z,t)=0)$ where $g_{2m}(z,t)\neq 0$,
$m\geq 1$ and $X_i$ is the weighted blow up 
$B_{(m,m,1,1)}X_{i+1}$.

\item ($cA_{>0}^+$, $\mult_0g$ odd) 
 Up to real analytic equivalence near $P$,
 $X_{i+1}\cong (x^2+y^2+g_{\geq 2m+1}(z,t)=0)$ where  $m\geq 1$,
$z^{2m+1}\in g$ and 
$z^it^j\not\in g$ for  $2i+j< 4m+2$. $X_i$ is the weighted blow up   
$B_{(2m+1,2m+1,2,1)}X_{i+1}$. 
\end{enumerate}
\end{enumerate}
\end{thm}

\begin{rem} The more precise results in sections 9--11 give a 
description of the various cases when $f_i$ is $\r$-small (though so
far I have not excluded some cases). 

The $\r$-trivial steps do not change anything in a neighborhood of the
real points, but it is in these steps that the full complexity of the
MMP appears. All the difficulties involving higher index terminal
singularities and flips are present, but they always appear in
conjugate pairs.

For the topological questions these have no effect, but in other
applications of (\ref{int.mmp.steps}) this should be taken into
account.
\end{rem}

\begin{rem}  The lists in (\ref{int.mmp.sings}) and
(\ref{int.mmp.steps}) are fairly short, but I do not see a simple
conceptual way of stating the results, let alone proving them by
general arguments.  The appearence of the singularities of type
$cE_6$ in (\ref{int.mmp.sings}) was rather unexpected for me.

The formulations also hide the cicumstance that there does not seem to
be a single method of excluding all other a priori possible cases. The
algebraic method of the proof of (\ref{int.mmp.sings}) ends with a
much longer list (\ref{ge.gwextr.thm}). The topological method excludes
many of these right away, but in a few cases several steps of the MMP
need to be analyzed.
\end{rem} 

\begin{say}[Method of the proof of (\ref{int.mmp.steps})]{\ }

The proof  relies on rather extensive computations.  The first step is
a classification of all 3--dimensional terminal singularities over
$\r$ and the study of their topological properties. This was carried
out in
\cite{rat1}. The next step is to gain  a good understanding of the
resolutions of these singularities. More precisely, we need to
understand the ``simplest" exceptional divisors in these resolutions.
(Simplicity is measured by the discrepancy, cf.\
(\ref{mmp.discr.def}).) Over $\c$ the first step in this direction is
\cite{Markushevich96}. A much more detailed study of such exceptional
divisors was completed by \cite{Hayakawa97}. Our main emphasis is over
$\r$, and it turns out that there is very little overlap between the
computations of \cite{Hayakawa97} and those in sections 9--11.
Nonetheless, the basic underlying principles are exactly the same.
\end{say}

\begin{ack}  I   thank M. Bestvina, S. Gersten, M. Kapovich and G.
Mikhalkin for answering my numerous questions about  3-manifold
topology and real algebraic geometry.  The existence of
$cE_6$ type points  in (\ref{int.mmp.sings}) was established with the
help of V. Alexeev.  I have received helpful comments and questions
from  A. Bertram, M. Fried, L. Katzarkov and B. Mazur.

Partial financial support was provided by  the NSF under grant number 
DMS-9622394. 
\end{ack}

\section{Applications and Speculations}

\subsection{Factorization of Birational Morphisms}

Let $f:Y\to X$ be a birational morphism between smooth and projective
varieties. It is a very old problem to factor $f$ as a composite of
``elementary" birational morphisms. In dimension 2 this is easy to do:
$f$ is the composite of blow ups of points. In dimension 3 and over
$\c$, the MMP factors $f$ as a composition of divisorial contractions
and flips, but these intermediate steps are  rather complicated and
not too well understood.

If $f:Y\to X$ is a birational morphism between smooth and projective
threefolds over $\r$, then one would like to get a factorization where
the intermediate steps are also defined over $\r$.   It turns out that
if $Y(\r)$ is orientable, the answer is very simple. As with minimal
models in general, the intermediate steps involve singular varieties
though in this case the real singularites are very mild.

\begin{defn} A real 3--fold $X$ is said to have a $cA_1$ singularity at
$0\in X(\r)$ if in  suitable real analytic cordinates $X$ can be given
by an equation $(\pm x^2\pm y^2\pm z^2\pm t^m=0)$ for a  suitable
choice of signs and   $m\geq 1$. 
\end{defn}

\begin{thm}\label{bir.morph.factor} Let $f:Y\to X$ be a birational
morphism between  smooth and projective threefolds over $\r$. Assume
that  
$Y(\r)$ satisfies the conditions (\ref{int.no.cond}). Then $f$ can be
factored as
$$
 f: Y=X_n\stackrel{f_n}{\to} X_{n-1}\to \cdots\to X_1 
\stackrel{f_1}{\to} X_0=X,
$$ where each $X_i$ has only $cA_1$ singularities at real points and 
the following is a complete list of possibilities for the $f_i$:
\begin{enumerate}

\item (smooth point blow up)  $f_i$ is the  blow up of a smooth point
$P\in X_{i-1}(\r)$.

\item (singular point blow up) $f_i$ is the blow up of a singular
point $P\in X_{i-1}(\r)$. 

\item (curve blow up) $f_i$ is the blow up of a real curve $C\subset
X_{i-1}$. $C$ has only finitely many real points, $X_{i-1}$ is smooth
at each of these and in suitable real analytic coordinates $C$ can be
written as $(z=x^2+y^{2m}=0)$.

\item ($\r$-trivial) $f_i$ is an isomorphism in a (Zariski)
neighborhood of the set of real points.
\end{enumerate}
\end{thm}

\begin{rem} As in (\ref{int.mmp.steps}),  it is  in  the $\r$-trivial
steps that the full complexity of the MMP appears. In particular,  
the $\r$-trivial steps  may be flips where the flipping curve has no
real points.
\end{rem}

Proof.  For purposes of induction we consider the more general case
when $X$ is allowed to have  $cA_1$-type singularities at real points
and terminal singularities at complex points. We assume  that $X$ is
$\q$-factorial (that is, a suitable multiple of every Weil divisor is
Cartier).  Run the real MMP for $Y$ over $X$ to obtain
$$
 f: Y=X_n\stackrel{f_n}{\map} X_{n-1}\map \cdots\map X_1 
\stackrel{f_1}{\to} X_0=X.
$$ The proof is by induction on the number of steps it takes the MMP
to reach $X$. 

The last step, $f_1:X_1\to X_0=X$,  is a contraction since we work over
$X$. 
 The possibilities for $f_1$  are described in  (\ref{int.mmp.steps}).
We are done  by induction  if $f_1$ is
$\r$-trivial or  a smooth point blow up. Assume that $f_1$ is a
singular point blow up.  Since $X_0$ has only $cA_1$ points, we are in
case (\ref{int.mmp.steps}.4a) with $m=1$. $f_1$ is the ordinary blow up
and by explicit computation we see that $X_1$ still has only 
$cA_1$ singularities.

The case when $f_1$ is $\r$-small (\ref{int.mmp.steps}.2) needs to be
studied in greater detail. $f_1$ can not be a g--extraction
(\ref{gw.g-e.def}) since $cA_1$ type points do not have g--extractions
other than the one listed above by (\ref{ge.gwextr.thm}). Thus $f_1$
is the blow up of a curve $C\subset X_0$. Moreover, 
$X_0$ is smooth along $C(\r)$ and 
$C$ is locally planar along $C(\r)$  by (\ref{sm.thm}). 
$C(\r)$ is finite  since $f_1$ is $\r$-small.

Pick any point $P\in C(\r)$ and assume that $C$ is given by real
analytic equations $(z=g(x,y)=0)$. By explicit computation, $B_CX_0$
has a unique singular point with equation
$(st-g(x,y)=0)$ which is equivalent to $(u^2-v^2-g(x,y)=0)$.

$X_1$ is an intermediate step of an MMP starting with $Y$, hence its
singularities are among those listed  (\ref{int.mmp.sings}).  Thus $g$
has multiplicity 2 and so it can be written as
$\pm x^2\pm y^r$.  Since $(g=0)$ has only the origin as its real
solution,   $g=\pm(x^2+y^{2m})$. \qed

\subsection{Application to the Nash Conjecture}

The main conclusion of (\ref{int.orient.thm}) and
(\ref{int.nonorient.thm}) is that if we want to understand the topology
of $X(\r)$ (say when it is orientable), it is sufficient to study the
topology of $X^*(\r)$ instead. $X^*$ has various useful properties,
depending on the conditions imposed on
$X$.

Consider, for instance, the original Nash question: what happens if $X$
is rational. Since the fifties  it has been understood that being
rational is a very subtle condition and it is very hard to work with.
\cite{KoMiMo92} introduced the much more general notion of being {\it
rationally connected}.  A $X$ is rationally connected if  two general
points of $X(\c)$ can be connected by an irreducible rational curve.
The lines show that $\p^n$ is rationally connected.

The structure theory of \cite{KoMiMo92} implies that a 3-fold  $X$ is
rationally connected iff
$X^*$ falls in one of 3 classes:
\begin{enumerate}
\item (Conic fibrations) There is a morphism (over $\r$) $g:X^*\to S$
onto a surface such that the general fiber is a conic. Correspondingly
there is a morphism $X^*(\r)\to S(\r)$ whose general fiber is $S^1$ or
empty. These cases will be studied in a subsequent paper.

\item (Del Pezzo fibrations) There is a morphism (over $\r$) $g:X^*\to
C$ onto a curve such that the general fiber is a Del Pezzo surface. If
$X(\r)$ is orientable, then this induces 
 a morphism $X^*(\r)\to C(\r)$ whose general fiber is a torus or a
union of some copies of $S^2$. These cases will be studied later.

\item (Fano varieties) The anticanonical bundle of $X^*$ is ample.
There is a complete list of such varieties if $X^*$ is also smooth
\cite{Iskovskikh80}.  Even if $X^*$ is known rather explicitly, a
topological description of
$X^*(\r)$ may not be easy. It would be interesting to work out at least
some of the cases, for instance hypersurfaces of degree 3 or 4 in
$\p^4$. (Mikhalkin pointed out that the degree 3 cases can be
understood using the classification of degree 4 real surfaces in
$\r\p^3$
\cite{Kharlamov76}.)

In general it is known that there are only finitely many families of
singular Fano varieties in dimension 3
\cite{Kawamata92}. Thus we can get only finitely many different
topological types for
$X^*(\r)$ in this case. 
\end{enumerate}

\subsection{Homology Spheres}

It is interesting to consider if we can get further simplifications of
the real MMP if we pose further restrictions on $X(\r)$. We may assume,
for instance, that $X(\r)$ is a homology sphere. This was in fact the
assumption I considered first. One can ask if under this assumption
$X(\r)\to X^*(\r)$ is a homeomorphism.

Unfortunately this is not the case. Consider for instance the singular
real threefold  $X^*$ given by affine equation
$$ x^2+y^2+z^2+(t-a_0)(t-a_m)\prod_{i=1}^{m-1}(t-a_i)^{2r}=0,
$$ where $a_0<a_1<\cdots<a_m$ are reals. This has $m-1$ singular points
of the form $x^2+y^2+z^2-u^{2r}=0$, which can be resolved by
$r$ successive blow ups. Resolving all singular points we obtain the
3--fold $X$. One can easily see that $X(\r)\sim S^3$, but
$\overline{X^*(\r)}$ is the disjoint union of $m$ copies of $S^3$. 

One may also study the types of singularities that occur if we pose
stronger restrictions on $X(\r)$.  It seems to me that the best one
can get is the following:

\begin{conj}\label{no.E6.conj}
 Let $X$ be a smooth projective 3--fold over $\r$. Assume that $X(\r)$
satisfies the conditions (\ref{int.no.cond}) and $X(\r)$ can not be
written as a connected sum
 with $S^1\times S^2$.

Let $X_i$ be any of the intermediate steps of the MMP over $\r$
starting with $X$ and $0\in X_i(\r)$ a real point. Then a neighborhood
of $0\in X_i$ is real analytically equivalent to one of the following
standard forms:
\begin{enumerate}

\item ($cA_0$) Smooth point.

\item ($cA_{>0}^+$) 
 $(x^2+y^2+g_{\geq 2}(z,t)=0)$, where  $\mult_0g$ is even and $g$ is
not everywhere negative in a punctured neighborhood of $0$.
\end{enumerate}
\end{conj}

In fact, most $cA_{>0}^+$-type singularities should  not occur.
  It is possible that one can write down a complete list. Also, one can
be more precise about how the singular points separate  $X_i(\r)$. 

The results in sections 8--11 come close to proving (\ref{no.E6.conj}),
but two points remain unresolved. In order to exclude $cE_6$ type
points,  one needs to show that the only possible g--extraction is the
one described in (\ref{cE6.g--extr.exist}). This should be a feasible
computation. The main problem is that in (\ref{cA+.multg-odd}) I could
not exclude certain $\r$-small contractions. I do not see how to deal
with this case.

\subsection{Beyond the Nash Conjecture}

One can refine the 3--dimensional Nash conjecture in two ways.
\medskip

First, one can study the topology of $X(\r)$ for other classes of real
algebraic varieties. The simplest cases may be those whose
 minimal models admit a natural fibration. This should be very helpful
in their topological study. One such class is elliptic threefolds,
where we have a morphism $X^*\to S$ whose general fiber is an elliptic
curve. A study of the singular fibers occurring in codimension 1 was
completed by \cite{Silhol84}.

Another, probably more difficult class are Calabi--Yau 3--folds. It
would be very interesting to find some connection between the topology
of $X(\r)$ and mirror symmetry.

The following question is consistent with the examples that I know:

\begin{question} Let $X$ be a smooth projective real 3--fold. Asume
that $X(\r)$ is hyperbolic. Does this imply that $X$ is of general
type?
\end{question}

\medskip

One can also start with a 3--manifold $M$ and look for a ``simple" real
projective 3--fold $X$ such that $X(\r)\sim M$. Ideally one would like
to find a solution where certain topological structures on $M$ are
reflected by the algebraic properties of $X$. 

There are hyperbolic 3--manifolds which embedd into
$\r^4$. Ths implies that they can be realized by real algebraic
hypersurfaces in $\r^4$. It would be interesting to find such examples.

\medskip

The methods of this paper require a very detailed study of the steps of
the MMP, which is currently feasible only in dimension 3.  It would
be, however, interesting to develop some examples in higher dimensions.

Example (\ref{connsum.exmp}) describing connected sum with $S^1\times
S^2$    should have interesting higher dimensional versions. There 
may be other, more complicated examples as well.

The first steps of the 4--dimensional MMP over $\c$ have been recently
classified by \cite{andwis}. It should be possible to obtain the
complete list over $\r$ and to study their topology.

\section{The Minimal Model Program over $\r$}

This section is intended to provide a    summary of the MMP over $\r$.
More generally, I   discuss the MMP over an arbitrary field $K$ of
characteristic zero, since there is no difference in the general
features. Conjecturally the whole program works in all dimensions but
at the moment it is only established in dimensions $\leq 3$.

\cite{koll87, koll90} provide general introductions. The minimal model
program for real algebraic surfaces is explained in detail in
\cite{ras}. 
 For more comprehensive treatments (mostly  over $\c$) see \cite{CKM88,
kolletal92, KoMo98}. 

One of the special features of the 3-dimensional MMP is that we have
to work with certain singular varieties in the course of the program.

\begin{defn}\label{mmp.qf.def}
 Let $X$ be a normal variety defined over a field $K$. A {\it (Weil)
divisor} over $K$ is a formal linear combination
$D:=\sum a_iD_i$  ($a_i\in \z$) of  codimension 1 subvarieties, each
defined and irreducible over $K$.  A {\it $\q$-divisor} is defined
similarly, except we allow $a_i\in \q$. A divisor $D$ is called {\it
Cartier} if it is locally definable by one equation and {\it
$\q$-Cartier} if  $mD$ is Cartier for some $m\in\n$. The smallest such
$m>0$ is called the {\it index} of $D$.

We say that $X$ is {\it factorial} (resp. {\it $\q$-factorial}) if
every Weil divisor is
 Cartier (resp.  $\q$-Cartier). 

A   divisor $D$ defined over $K$ is Cartier (resp.  $\q$-Cartier) iff
it is Cartier (resp.  $\q$-Cartier)  after some field extension.
However, a variety may be $\q$-factorial over $K$ and not
$\q$-factorial over $\bar K$. For instance, the cone
$x^2+y^2+z^2-t^2$ is factorial over $\r$ but not over $\c$.  (For
instance, $(x-\sqrt{-1}y=z-t=0)$ is a not $\q$-Cartier.)
\end{defn}

\begin{defn}\label{mmp.KX.def}
 For a normal variety $X$, let
 $K_X$ denote its {\it canonical class}. $K_X$ is a linear equivalence
class of Weil divisors. The corresponding reflexive sheaf $\o_X(K_X)$
is isomorphic to the {\it dualizing sheaf}
$\omega_X$ of
$X$. 

The index of $K_X$ is called the {\it index } of $X$. 
\end{defn}

\begin{defn}\label{mmp.discr.def} Let $X,Y$ be normal varieties and 
$f:Y\to X$  a birational morphism with exceptional set $\ex(f)$.  Let
$E_i\subset \ex(f)$ be the exceptional divisors. If
$mK_X$ is Cartier, then 
$f^*\o_X(mK_X)$ is defined and there is a natural isomorphism
$$ f^*\o_X(mK_X)|(Y\setminus \ex(f))\cong 
\o_Y(mK_Y)|(Y\setminus \ex(f)).
$$
 Hence there are integers $b_i$ such that
$$
\o_Y(mK_Y)\cong f^*\o_X(mK_X)(\sum b_iE_i).
$$ Formally divide by $m$ and write this as
$$ K_Y\equiv f^*(K_X)+\sum a(E_i, X)E_i,\qtq{where $a(E_i,X)\in \q$.}
$$ The rational number $a(E_i,X)$ is called the {\it discrepancy} of
$E_i$ with respect to $X$.

The closure of $f(E_i)\subset X$ is called the {\it center} of $E_i$ on
$X$. It is denoted by $\cent_XE_i$.

If $f':Y'\to X$ is another birational morphism and $E'_i:=(f'\circ
f^{-1})(E_i)\subset Y'$ is a divisor then $a(E'_i,X)=a(E_i,X)$ and
$\cent_XE_i=\cent_XE'_i$.  Thus the discrepancy and the center depend
only on the divisor up to birational equaivalence, but not on the
particular variety where
 the divisor appears.
\end{defn}

\begin{defn} Let $X$ be a normal variety such that $K_X$ is
$\q$-Cartier.  We say that $X$ is {\it terminal}  (or that it has {\it
terminal singularities}) if for every $f:Y\to X$, the discrepancy of
every exceptional divisor is positive.
\end{defn}

The following result makes it feasible to decide if  $X$ is terminal
or not. 

\begin{lem}\label{mmp.term.lem} For a normal variety $X$ the following
are equivalent:
\begin{enumerate}
\item $X$ is terminal,
\item $a(E,X)>0$ for every   resolution of singularities $f:Y\to X$
and for every  exceptional divisor $E\subset \ex(f)$.
\item  There is a   resolution of singularities
$f:Y\to X$  such that  $a(E,X)>0$  for every  exceptional divisor
$E\subset \ex(f)$.\qed
\end{enumerate}
\end{lem}

\begin{exmp}\label{mmp.comp.discr.exmp}
 It is frequently not too hard to compute discrepancies. Assume for
instance that $X$ is a hypersurface defined by
$(F(x_1,\dots,x_n)=0)$. A local generator of $\o_X(K_X)$ is given by
any of the forms
$$
\eta_i:=\frac1{\partial F/\partial x_i}dx_1\wedge\cdots\wedge
dx_{i-1}\wedge dx_{i+1}\wedge\cdots\wedge dx_n.
$$ Let $f:Y\to X$ be a resolution of singularities and $P\in Y$ a point
with local coordinates $y_1,\dots,y_{n-1}$. $f$ is given by coordinate
functions
$x_i=f_i(y_1,\dots,y_{n-1})$ and so we can write
\begin{eqnarray*} f^*\eta_n&=&f^*\left(\frac1{\partial F/\partial
x_n}\right)
\operatorname{Jac} dy_1\wedge\cdots\wedge dy_{n-1},
\qtq{where}\\
\operatorname{Jac}&=& \operatorname{Jac}
\left(\frac{f_1,\dots,f_{n-1}}{x_1,\dots,x_{n-1}}\right)
\end{eqnarray*} denotes  the determinant of the Jacobian matrix. 
Hence the discrepancies can be computed as the order of vanishing of
the Jacobian minus the order of vanishing of
$f^*(\partial F/\partial x_n)$. 

 If $X$ is smooth then we conclude  that $a(E,X)\geq 1$ for every 
exceptional divisor. Thus smooth varieties are terminal.
\end{exmp}

Next we define various birational maps which have special role in the
MMP.

\begin{defn}\label{mmp.extremal.def}
 Let $X$ be a variety over $K$ and assume that $K_X$ is
$\q$-Cartier. A proper morphism $g:X\to Y$ is called an {\it extremal
contraction} if the following  conditions hold:
\begin{enumerate}
\item $g_*\o_X=\o_Y$,
\item $X$ is $\q$-factorial,
\item Let $C\subset X$  be any irreducible curve  such that
$g(C)=\mbox{point}$. Then a 
 $\q$-divisor $D$ on $X$ is the pull back of a $\q$-Cartier
$\q$-divisor  $D'$ on $Y$  iff $(D\cdot C)=0$. (Necessarily,
$D'=g_*(D)$.)
\end{enumerate}
\end{defn}

\begin{defn}\label{mmp.excontypes.def}
 Let $g:X\to Y$ be an extremal contraction. 

We say that $g$ is of {\it fiber type} if $\dim Y<\dim X$.

 We say that $g$ is a {\it divisorial} contraction if the exceptional
set $\ex(g)$ is the support of $\q$-Cartier divisor. In this case
$\ex(g)$ is irreducible over $K$.

 We say that $g$ is a {\it small} contraction if
 $\dim \ex(g)\leq \dim X-2$.

One can see that  every extremal contraction is in one of these 3
groups.
\end{defn}

\begin{defn}\label{mmp.KXneg.def}
 A proper morphism $f:X\to Y$ is called {\it $K_X$-negative} if $-K_X$
is $f$-ample.
\end{defn}

\begin{defn}\label{mmp.flip.def}  Let 
$f:  X  \to   Y$  be a small $K_X$-negative extremal contraction.
 A variety $X^+$ together with a proper birational morphism 
$f^+: X^+
\to  Y$  
	is called a {\it flip} of  $f$  if
\begin{enumerate} 
\item  $K_{X^+}$  is $\q$-Cartier, 
\item  $K_{X^+}$ is $f^+$-ample, and
\item  the exceptional set  $\ex(f^+)$ has codimension 	at least two
in $X^+$.  
\end{enumerate} By a slight   abuse of terminology, the rational map
$\phi: X \map X^+$           is also
 called a flip.  A flip gives the following diagram:
$$
\begin{array}{rcl} X &\stackrel{\phi}{\map} &X^+\\
\mbox{$-K_{X}$ is $f$-ample} &\searrow \quad\swarrow &
\mbox{$K_{X^+}$ is $f^+$-ample}\\ &Y&
\end{array}
$$ It is not hard to  see that a flip is unique and the main question
is its existence.
\end{defn}

We are ready to state the 3-dimensional MMP over an arbitrary field:

\begin{thm}[MMP over $K$]\label{mmp.mmp.thm}
  Let $X$ be a smooth projective 3-fold defined over a field $K$ (of
characteristic zero). Then there is a sequence
$$
 X=X_0\stackrel{f_0}{\map} X_1\map \cdots \map X_i
\stackrel{f_i}{\map} X_{i+1}\map \cdots \stackrel{f_{n-1}}{\map}
X_n=:X^*
$$ with the following properties
\begin{enumerate}
\item  Each $X_i$  is a terminal projective 3-fold over $K$ which is
$\q$-factorial over $K$.
\item Each $f_i$ is either a $K_X$-negative divisorial  extremal
contraction or the  flip of a $K_X$-negative small  extremal
contraction.
\item One of the following holds for $X^*$:
\begin{enumerate}
\item
 $K_{X^*}$ is nef (that is $(C\cdot K_{X^*})\geq 0$ for any curve
$C\subset X^*$), or
\item  there is a fiber type extremal contraction $X^*\to Z$. 
\end{enumerate}
\end{enumerate}
\end{thm}

\begin{rem} For the purposes of this paper one can handle the  MMP as
a black box. It is sufficient to know that it works, but I will  use
very few of its finer properties. In particular, there is no need to
know anything about flips beyond  believing their existence.

The rest of the section is devoted to explicitly stating all further
results from minimal model theory  that I use later. The most
significant among these is the  classification of terminal 3-fold
singularities over nonclosed fields, established in \cite{rat1}. 
\end{rem}

\begin{notation}\label{mmp.ps.not}
 For a field $K$ let
$K[[x_1,\dots,x_n]]$ denote the ring of formal power series in $n$
variables over
$K$. For $K=\r$ or
$K=\c$, let $K\{x_1,\dots,x_n\}$ denote the ring of those formal power
series which converge in some  neighborhood of the origin.

For a power series $F$, $F_d$ denotes the degree $d$ homogeneous part.
The multiplicity, denoted by $\mult_0F$, is the smallest $d$ such that
$F_d\neq 0$. If we write a power series as $F_{\geq d}$ then it is
assumed that its multiplicity is at least $d$.

 For $F\in \r\{x_1,\dots,x_n\}$ let $(F=0)$ denote the germ of its
zero set in $\c^n$ with its natural real structure.  I always think of
it as a complex analytic germ with a real structure and not just as a
real analytic germ in $\r^n$. 

$(F=0)/\frac1{n}(a,b,c,d)$ means the following. Define a
$\z_n$-grading of $\c\{x,y,z,t\}$ by $x\mapsto a, y\mapsto b, z\mapsto
c,t\mapsto d$.  If $F$ is graded homogeneous, then
$(F=0)/\frac1{n}(a,b,c,d)$ denotes the germ whose ring of holomorphic
functions is the ring of grade zero elements of
$\c\{x,y,z,t\}/(F)$. 

If $(F=0)$ is terminal then $n$ coincides with the index 
(\ref{mmp.KX.def})   of the singularity.
\end{notation}

\begin{exmp} In case $X=(x^2+y^2+z^2+ t^2=0)/\frac12(1,1,1,0)$
 the ring is
$$
\o_X=\c\{x^2,y^2,z^2,t,xy,yz,zx\}/(x^2+y^2+z^2+ t^2),
$$ with the natural real structure. 

$X$ can also be realized as the image of  the hypersurface
$(x^2+y^2+z^2+ t^2=0)$ under the map
$$
\phi: \c^4\to \c^7:\quad (x,y,z,t)\mapsto (x^2,y^2,z^2,t,xy,yz,zx),
$$ which has degree 2 over its image.

Although
$(x^2+y^2+z^2+ t^2=0)$ has only the origin as its real solution,
$X$ has plenty of real points. Indeed, any real solution of
$x^2+y^2+z^2- t^2=0$ gives a {\it real} point 
$P=\phi(\sqrt{-1}x,\sqrt{-1}y, \sqrt{-1}z, t)\in X(\r)$. 
$\phi^{-1}(P)$ is a pair of conjugate points on the hypersurface 
$(x^2+y^2+z^2+ t^2=0)$. All the real elements of
$\o_X$ take up real values at  $P$.

This way we see that $X(\r)$ is a cone over 2 copies of $\r\p^2$. 
\end{exmp}

The following is a summary of the classification of terminal
singularities obtained in \cite{rat1}.  As it turns out, the
classification closely follows the earlier results over algebraically
closed fields. The choice of the subdivison into cases is dictated by
the needs of the proof in sections 9--11, rather than the internal
logic of the classification.

\begin{thm}\label{mmp.ts.thm}  Let $X$ be a real algebraic or analytic
3-fold  and
$0\in X(\r)$ a real point.  Then $X$ has a terminal singularity at $0$
iff a neighborhood of
$0\in X$ is real analytically equivalent to one of the following:
$$
\begin{tabular}{ll} name & \qquad equation \\
 $cA_0$     &$(t=0)$\\
 $cA_1$      &$(x^2+y^2\pm z^2\pm t^m=0)$\\
 $cA_{>1}^+$     &$(x^2+y^2+g_{\geq 3}(z,t)=0)$  \\
 $cA_{>1}^-$     &$(x^2-y^2+g_{\geq 3}(z,t)=0)$ \\
 $cD_4$     &$(x^2+f_{\geq 3}(y,z,t)=0)$, where $f_3\neq l_1^2l_2$ for
linear forms $l_i$\\
 $cD_{>4}$     &$(x^2+y^2z+f_{\geq 4}(y,z,t)=0)$,\\
 $cE_6$     &$(x^2+y^3+yg_{\geq 3}(z,t)+h_{\geq 4}(z,t)=0)$, where
$h_4\neq 0$\\
 $cE_7$     &$(x^2+y^3+yg_{\geq 3}(z,t)+h_{\geq 5}(z,t)=0)$, where
$g_3\neq 0$\\
 $cE_8$     &$(x^2+y^3+yg_{\geq 4}(z,t)+h_{\geq 5}(z,t)=0)$, where
$h_5\neq 0$\\
 $cA_0/n$     & $(t=0)/\frac1{n}(r,-r,1,0)$ where $n\geq 2$ and
$(n,r)=1$\\
 $cA_1/2$     & $(x^2+y^2\pm z^n\pm t^m=0)/\frac12(1,1,1,0)$ where
$\min\{n,m\}=2$\\
 $cA_{>1}^+/2$     & $(x^2+y^2+f_{\geq 3} (z,t)=0)/\frac12(1,1,1,0)$\\
 $cA_{>1}^-/2$     & $(x^2-y^2+f_{\geq 3}(z,t)=0)/
\frac12(1,1,1,0)$\\
 $cA/n$     & $(xy+f(z,t)=0)/\frac1{n}(r,-r,1,0)$ where $n\geq 3$ and
$(n,r)=1$\\
 $cAx/2$     &  $(x^2\pm y^2 +f_{\geq 4} (z,t)=0)/\frac12(0,1,1,1)$\\
 $cAx/4$     & $(x^2\pm y^2 +f_{\geq 2}(z,t)=0)/\frac14(1,3,1,2)$
where $f_2(0,1)=0$\\
 $cD/2$     &$(x^2+f_{\geq 3}(y,z,t)=0)/\frac12(1,1,0,1)$\\
 $cD/3$     &$(x^2+f_{\geq 3}(y,z,t)=0)/\frac13(0,1,1,2)$ where
$f_3(0,0,1)\neq 0$\\
 $cE/2$     &$(x^2+y^3+f_{\geq 4}(y,z,t)=0)/\frac12(1,0,1,1)$\\
\end{tabular}
$$
\end{thm}

\section{The Topology of Real Points and the MMP}

 Starting with a  projective variety 
$X$ over $\r$, let us run the MMP over $\r$. We obtain a sequence of
birational maps
$$ X=X_0\map X_1\map \cdots\map X_i\stackrel{f_i}{\map}
X_{i+1}\map\cdots\map X^*.
$$ These in turn induce (not necessarily everywhere defined) maps
between the sets of real points
$$ X(\r)=X_0(\r)\map  \cdots\map X_i(\r)\stackrel{f_i}{\map}
X_{i+1}(\r)\map\cdots\map X^*(\r).
$$ Our aim is to see if there is a way of describing $X(\r)$ in terms
of 
$X^*(\r)$ and  a local description of the maps 
$X_i(\r)\map X_{i+1}(\r)$ in a neighborhood of their exceptional sets.

\begin{prop}\label{mmpt.cases.prop}
 Every step $f_i$ of the MMP over
$\r$ is among the following five:
\begin{enumerate}
\item (divisor--to--point) $f_i$ contracts a geometrically irreducible
divisor $E_i\subset X_i$ to a point $P_{i+1}\in X_{i+1}(\r)$.

\item (divisor--to--curve) $f_i$ contracts a geometrically irreducible
divisor $E_i\subset X_i$ to a real curve $C_{i+1}\subset X_{i+1}$.

\item ($\r$-small) $f_i:X_i(\r)\to X_{i+1}(\r)$ collapses a 1-complex
to   points and is a homeomorphism elsewhere.

\item (flip) $f_i$ is the flip of a curve $C_i\subset X_i$.

\item ($\r$-trivial) $f_i$ is an isomorphism in a (Zariski)
neighborhood of the set of real points.
\end{enumerate}
\end{prop}

Proof. If $f_i$ is a flip then we have case (4). Thus we may assume
that $f_i$ is the contraction of a divisor
$E_i\subset X_i$ and $E_i$ is irreducible over $\r$. If $E_i$ is
irreducible over $\c$ then we have one of the cases (1--2). If $E_i$
is reducible over $\c$ then $E_i(\r)$ is a 1-complex by
(\ref{mmpt.red.small.lem})  and so we are in case (3).

Any of the above cases can also be of type (5).\qed

\begin{lem}\label{mmpt.red.small.lem}
 Let $X$ be an n-dimensional scheme over $\r$ (that is, an algebraic
variety possibly  with several irreducible components and with
singularities).  Assume that if $X_i\subset X$ is any $\r$-irreducible
component then
$X_i$ is reducible over $\c$. Then
$X(\r)= (\sing X)(\r)$, that is, every real point is singular. In
particular, $\dim X(\r)\leq n-1$. 
\end{lem}

Proof.  Assume that $P\in X(\r)$ is a smooth real point. Then $P$ lies
on    a unique irreducible component $Y\subset X_{\c}$, thus $Y$ is
invariant under complex conjugation. So
$Y$ is an  irreducible real component which stays irreducible over
$\c$, a contradiction.\qed 
\medskip

Each of the 5  steps $X_i(\r)\map X_{i+1}(\r)$  have different
topological behaviour. The following informal  discussion   intends to
emphasize their main features.

\begin{say}[Divisor--to--point]   Let $M=X_i(\r)$ be a 3-complex 
(with only finitely many singular points) and  
$F=E_i(\r)\subset M$  a 2-complex.  We collapse $F$ to a point:
$$
\begin{array}{ccc} F& \subset & M\\
\downarrow && \ \ \downarrow f\\ P&\in &N.
\end{array}
$$ In practice we are frequently able to describe a regular
neighborhood
$F\subset U\subset M$  (this is a local datum) and by assumption we
know a regular neighborhood $P\in V\subset N$. Thus we see that M is
obtained from $U$ and $N\setminus
\inter V$ by gluing them together along the boundaries
$\partial U$ and $\partial V$.

The gluing is determined by a PL-homeomorphism $\phi:
\partial U\to \partial V$. Thus, besides knowing
$U$ and
$N$, we also need to know $\phi$  up to PL-isotopy.  If one of the
connected components of $\partial U$   has genus at least 2, this is a
very hard problem. In fact, as the example of Heegard splittings shows
(cf.
\cite[Ch.2]{Hempel76}), the choice of $\phi$ is  usually the most
significant information.   Unfortunately, $\phi$  can be described
only in terms of  global data.

If $\partial U$ is a union of $m$ copies of $S^2$, then
$\phi$ is classified by an element of the symmetric group on
$m$ elements (which $S^2$ maps where) and a sign for each
$S^2$ (describing whether the map is orientation preserving or
reversing on that $S^2$). Hence, knowing $U$ and  $N$, we can
determine $M$ up to  finite ambiguity.

In many cases $U$ is so simple that  different choices of
$\phi$  give the same $M$, giving even fewer possibilities for $M$.

If $P\in N$ is an isolated singular point, then $\partial V$ is a
union of spheres iff $\bar N$ (the topological normalization of $N$)
is a manifold.

The situation is similarly simple if $\partial U$  is a union of 
copies of $\r\p^2$ and of $S^2$, and still  manageable if $\partial
U$  also contains  tori and Klein bottles. For us these more general
cases do not come up.
\end{say}

\begin{say}[Divisor--to--curve] This  time we construct  it bottom up.
Assume for simplicity  that 
$N=X_{i+1}(\r)$ is a 3-manifold   and  
$L=C_{i+1}(\r)\subset N$  a link.   The projectivized normal bundle is
an $S^1$-bundle  $S\to L$. The blow up of $L$ in $N$  replaces $L$ by
$S$ to obtain:
$$
\begin{array}{ccc} S& \subset & M\\
\downarrow && \ \ \downarrow f\\ L&\in &N.
\end{array}
$$ (In general $N$ may have finitely many singular points and $L$ is
only a 1-complex, but I believe that  a similar description is
possible in all cases.) 

Here $M$ is uniquely  determined, once we know $N$ and $L$. By
assumption we know $N$ but $L$ is a free choice.   The
Jaco--Johannson--Shalen decomposition (cf. \cite[p.483]{Scott83})
shows that in most cases $B_LN$ determines $M\setminus L$. Thus the
description of all possible $B_LN$ is essentially equivalent to the
description of all links.

For us
$L$ has to come from an algebraic curve, thus we are led to the
question: Which links in a real algebraic 3-fold can be realized by
algebraic curves? In some cases  every link is realized (cf.
\cite{AK81}), thus we again run into  a hard topological problem.

So  $M$ can be  described in terms of $N$, though the answer depends
on the choice of a link, which is a very complicated object.
\end{say}

\begin{say}[$\r$-small contraction]    $N=X_{i+1}(\r)$ is obtained
from $M=X_i(\r)$ by collapsing a 1-complex
$C=(\sing E_i)(\r)=E_i(\r)$  to a point:
$$
\begin{array}{ccc} C& \subset & M\\
\downarrow && \ \ \downarrow f\\ P&\in &N.
\end{array}
$$ If the normalizations  $\bar M$ and $\bar N$ are manifolds, then 
we see in (\ref{top.alg.homeo.lem}) that a suitable small perturbation
of
$f$ is a  homeomorphism between $\bar M$ and $\bar N$.  Thus this
step  (which is actually more complicated from the point of view of
algebraic geometry than the previous two cases)  is easy to analyze
topologically.
\end{say}

\begin{say}[Flip]  Assume for simplicity that 
 $M=X_i(\r)$ is an orientable  3-manifold   and  
$C(\r)\sim S^1$.   $N=X_{i+1}(\r)$ is obtained from $M$ by a surgery
along $S^1$.  The boundary of a regular neighborhood of $S^1$ is
$S^1\times S^1$, and the surgery is determined by a diffeomorphism of
$S^1\times S^1$ up to isotopy. These are classified by $SL(2,\z)$. (In
general $M$ may have finitely many singular points and $C(\r)$ is  a
1-complex, but I believe that  a similar description is possible in
all cases.) 

A complete classification of flips is known \cite{KoMo92}, thus it
should be possible to compute the resulting  diffeomorphism of
$S^1\times S^1$.

Here again we run into a global problem. $S^1\subset M$ may be very
knotted, and the result of the surgery depends mostly on the knot
$S^1\subset M$. The usual descriptions of flips characterize a complex
analytic neighborhood of $C$, thus they say nothing about how its real
part is knotted. From the point of view of algebraic geometry, this 
is a global invariant. 

We have the additional problem that  flipping curves are  rigid
objects,  thus we can not hope to get a  flipping curve by
approximating  a real curve algebraically. Furthermore,   it is very
hard to determine which curves are obtained by a flip. (Even if $Z$ is
a smooth complex 3-fold and $C\subset Z$ a smooth
$\c\p^1$, I know of no practical way of determining if
$C\subset Z$ is obtained as  a result of a flip.)
\end{say}

\begin{say}[Conclusion]\label{mmpt.concl} 
 Start with a  projective 3-fold 
$X$ over $\r$ and run the MMP over $\r$:
$$ X=X_0\map X_1\map \cdots\map X_i\stackrel{f_i}{\map}
X_{i+1}\map\cdots\map X^*.
$$ If we would like to understand the topology of
$X(\r)$ in terms of  $X^*(\r)$, then we have to ensure that the MMP
has the following properties:
\begin{enumerate}
\item 
$\overline{X_i(\r)}$ is a manifold for every $i$.
\item Each $f_i$  is either $\r$-trivial or $\r$-small or a
divisor--to--point contraction.
\end{enumerate}
\end{say}

(\ref{int.mmp.steps}) asserts that both of these conditions can be
satisfied by imposing certain mild conditions on the topology of
$X(\r)$.

\section{The Topology of Divisorial Contractions}

The aim of this section is to describe some examples where the change
of the topology  of a real algebraic variety under a divisorial
contraction  can be readily understood by topological methods.

\begin{notation}\label{5.notation}
 The disjoint union of two topological spaces is denoted by $M\uplus
N$.  Direct product is denoted by $M\times N$. The unique nontrivial
$S^2$-bundle over $S^1$ is denoted by 
$S^1\tilde{\times} S^2$. This is obtained from $[0,1]\times S^2$ by
indentifying the 2 ends via an orientation reversing homeomorphism.

 Homeomorphism of two topological spaces is denoted by $M\sim N$. 
\end{notation}

We start with the study of $\r$-small contractions:

\begin{lem}\label{top.homeo.lem}
 Let $f:M\to N$ be a   proper PL-map between PL-manifolds of dimension
$n\geq 3$. Assume that there is a 1-complex
$C\subset M$ and a finite set of points $P\subset N$ such that
$f:M\setminus C\to N\setminus P$ is a PL-homeomorphism.

Then $M$ and $N$ are PL-homeomorphic (by a small perturbation of $f$).
\end{lem}

Proof. If $C$ is collapsible to points, then  a regular neighborhood
of $C$ in $M$ is a union of disjoint  $n$-cells
\cite[1.8]{Hempel76} and we are done. 

In order to see that $C$ is collapsible to points, we may assume that
$P$ is a point and
$N=S^n$. Thus $M$ is also a compact PL-manifold.
$M$ is orientable outside the codimenison $\geq 2$ subset 
$C$, hence it is orientable.  Consider the exact homology sequences
$$
\begin{array}{ccccccc} H_i(C)&\to &H_i(M)&\to & H_i(M,C)&\to &
H_{i-1}(C)\\
\downarrow &&\downarrow &&\downarrow &&\downarrow \\ H_i(P)&\to
&H_i(S^n)&\to & H_i(S^n,P)&\to & H_{i-1}(P)
\end{array}
$$ We compute $H_i(M,C)=H_i(S^n,P)$ from the second sequence and
substitute into the first to obtain that
$$ H_1(C)\cong H_1(M), \qtq{and} 0=H_{n-1}(C)\cong H_{n-1}(M).
$$ By Poincar\'e duality we conclude that $H_1(C)=0$, thus
$C$ is contractible.
\qed

\begin{cor}\label{top.alg.homeo.lem}
 Let $f:X\to Y$ be a morphism a $n$-dimensional real algebraic
varieties,
$n\geq 3$. Assume that 
\begin{enumerate}
\item $\overline{X(\r)}=M\uplus R$ and $\overline{Y(\r)}=N\uplus R'$
where $M,N$  are PL-manifolds and $\dim R,\dim R'<n$.
\item $f$ induces an isomorphism $R\cong R'$.
\item  $\ex(f)(\r)$ is a 1-complex.
\end{enumerate}

\noindent Then $\overline{X(\r)}$ is PL-homeomorphic to
$\overline{Y(\r)}$.
\end{cor}

Proof. Set $C=\ex(f)(\r)$ and $\bar C\subset M\subset 
\overline{X(\r)}$ its preimage. Since $\bar C$ has dimension 1, there
is a one--to--one correspondence between the connected components of
$\bar C$ and the connected components of the boundary of a regular
neighborhood of $\bar C$. Hence $f$ lifts to a morphism $\bar f:
\overline{X(\r)}\to
\overline{Y(\r)}$.
 Thus (\ref{top.homeo.lem}) implies (\ref{top.alg.homeo.lem}), and the
homeomorphism is given  by a small perturbation of
$\bar f$.\qed

\begin{rem} A real algebraic curve is a union of copies of
$S^1$.  The proof of (\ref{top.homeo.lem}) shows that the preimage of
$\ex(f)(\r)$ in $\overline{X(\r)}$ is contractible. Thus if
$f$ itself is not a homeomorphism, then $X(\r)$  is not a manifold.
\end{rem}

Next we look at divisor--to--point contractions.

\begin{prop}\label{top.normsurf.nbd.boundary}
 Let $M$ be a 3-dimensional PL-manifold and $F\subset M$ a connected 
2-complex with only finitely many singular points.  Let
$F\subset \inter U\subset M$ be a regular neighborhood of
$F$. Then
$$
\dim H_1(\bar F,\q)\leq \dim H_1(\partial U,\q),
$$ and strict inequality holds unless every connected component of
$\bar F$ is one of the following:
\begin{enumerate}
\item $S^2$ or $\r\p^2$,
\item a one-sided $S^1\times S^1$,
\item a one-sided Klein bottle whose neighborhood is not orientable.
\end{enumerate}
\end{prop}

Proof.  Let $F$ be a compact 2-complex  with only finitely many
singular points. Its normalization  $\bar F$ can be written as
$F^{(2)}\uplus F^{(1)}$ where  $F^{(2)}$ is a compact 2-manifold and a
$ F^{(1)}$ is a 1-complex.

 Pick a point $P\in F$ whose link in $F$ consists of at least 2
circles.  Locally $F$ looks like the cone over  parallel plane
sections  $(z=a_i)\cap (x^2+y^2+z^2=1)$ of the unit sphere in
$\r^3$ (plus a few 1-cells). By a homotopy we can replace this by the
parallel plane sections of the unit ball
$(z=a_i)\cap (x^2+y^2+z^2\leq 1)$
 and add the interval $[\min_i\{a_i\},\max_i\{a_i\}]$ on the
$z$-axis. This does not change the boundary of the regular
neighborhood. Thus we may assume that $F^{(2)}\to M$ is an embedding.

Let us take a point or a 1-cell $e$ in $F^{(1)}$.  If  $e$ does not
intersect the rest  of $F$, then a regular neighborhood of
$e$ is a 3-cell.  $e$ can be deleted from $F$ without changing the
inequality.

If  $e$ intersects the rest  of $F$ in one endpoint only, then we can
delete $e$ from $F$ without changing the regular neighborhood. 

If $e$ intersects the rest  of $F$ at both endpoints, then removing
$e$ creates a new 2-complex $F'$, and $F^{(2)}={F'}^{(2)}$. Let
$F'\subset\inter U'$ be its regular neighborhood. 
$\partial U$ is obtained from $\partial U'$ by attaching a handle
$[0,1]\times S^1$. Thus $H_1(\partial U)\geq H_1(\partial U')$, and it
is sufficient to verify our inequality for $F'$.

At the end we are reduced to the situation when $F$ is  the disjoint
union of embedded 2-manifolds, and it is sufficient to check the
ineqality for each connected component of $F$ separately.
$\partial U \to F$ is a 2 sheeted cover, thus
$H_1(\partial U)\geq H_1(F)$ with equality only if
$F\sim S^2, F\sim \r\p^2, F\sim S^1\times S^1\sim \partial U$  or
$F$ and $\partial U$ are both Klein bottles.\qed

\begin{prop}\label{top.normsurf.collapse}
 Let $M$ be a 3-dimensional PL-manifold and $F\subset M$ a compact 
2-complex with only finitely many singular points.  Let $0\in N$ be
obtained from $M$ by collapsing $F$ to a point. Assume that
$\bar N$ is a 3-manifold. Then
 $M$ can be obtained from $\bar N$ by repeated application of the
following operations:
\begin{enumerate}
\item taking connected sums of connected components,
\item taking connected sum with $S^1\times S^2$,
\item taking connected sum with $S^1\tilde{\times} S^2$, or
\item taking connected sum with $\r\p^3$.
\end{enumerate}
\end{prop}

Proof.  We use the notation of (\ref{top.normsurf.nbd.boundary}) and
of its proof. Let
$F\subset \inter U\subset M$ and $0\in \inter V\subset N$ be regular
neighborhoods such that $U=f^{-1}(V)$.  Then $\partial U=\partial V$.
Since
$\bar N$ is a manifold, this implies that $\partial U$ is a union of
2-spheres. We also see that $\bar N$ is obtained from
$M\setminus\inter U$ by attaching a 3-ball to each 
$S^2$ in $\partial U$. 

As in the proof of (\ref{top.normsurf.nbd.boundary}) we may assume that
$F^{(2)}\to M$ is an embedding.

If $e$ is a point or a 1-cell  in $F^{(1)}$ which   intersects the
rest  of
$F$ in zero or one point only, then we can delete $e$ from $F$.

If $e$ intersects the rest  of $F$ at both endpoints, then removing
$e$ creates a new 2-complex $F'$ such that $\bar F=\bar F'$. Let
$F'\subset\inter U'$ be its regular neighborhood. 
$\partial U$ is obtained from $\partial U'$ by attaching a handle
$[0,1]\times S^1$.  The two ends $\{0\}\times S^1$ and 
$\{1\}\times S^1$ can not attach to the same connected component of 
$\partial U'$ since that would create a torus or a Klein bottle in
$\partial U$. Thus $\partial U'$ has one more  copies of
$S^2$ than 
$\partial U$. 

  $\bar N$  is obtained from  $\bar N'$  by collapsing the image of
$e$ to a point, hence 
$\bar N$  and   $\bar N'$  are homeomorphic by (\ref{top.homeo.lem}). 

At the end we are reduced to the situation when $F$ is  the disjoint
union of embedded copies of $S^2$ and
$\r\p^2$.  An  $S^2$  is necessarily 2-sided. Removing it from $F$
corresponds to taking connected sums of connected components (if $S^2$
separates $M$) or to  taking connected sum with $S^1\times S^2$ or 
$S^1\tilde{\times} S^2$ (if $S^2$ does not separate
$M$) (cf. \cite[Chap.\ 3]{Hempel76}).

If $\r\p^2$ is 2-sided, then the boundary of its  regular neighborhood
consists of two copies of $\r\p^2$, so this can not happen. A 1-sided
$\r\p^2$ corresponds to taking connected sum with
$\r\p^3$. \qed

\begin{cor}\label{top.alg.normsurfup.lem}
 Let $f:X\to Y$ be a morphism a real algebraic $3$-folds. Assume that 
\begin{enumerate}
\item $\overline{X(\r)}$ and $\overline{Y(\r)}$ are PL-manifolds, and
\item  $\ex(f)$ is a  geometrically irreducible normal surface which is
contracted to a point. 
\end{enumerate}

\noindent Then $\overline{X(\r)}$  can be obtained from
$\overline{Y(\r)}$ by repeated application of the following
operations: 
\begin{enumerate}
\setcounter{enumi}{2}
\item removing an isolated point from $\overline{Y(\r)}$,
\item taking connected sums of connected components,
\item taking connected sum with $S^1\times S^2$,
\item taking connected sum with $S^1\tilde{\times} S^2$, or
\item taking connected sum with $\r\p^3$.
\end{enumerate}
\end{cor}

Proof.  If $\ex(f)(\r)=\emptyset$ then  the image of $\ex(f)$  is an
isolated real point of $Y(\r)$ which has to be thrown away to obtain 
$X(\r)$.  If $\ex(f)(\r)\neq\emptyset$, then isolated points of
$X(\r)$ correspond to isolated points of $Y(\r)$, hence they can be
ignored.

Let $M$ be the topological normalization of
$X(\r)\setminus(\mbox{isolated points})$,
$N$ the topological normalization of $Y(\r)\setminus(\mbox{isolated
points})$ and  
$F$   the preimage of $\ex(f)(\r)$ in $M$.  $F$ is a 2-complex with
isolated singularities since $\ex(f)$ is normal.

Thus  (\ref{top.alg.normsurfup.lem}) follows from 
(\ref{top.normsurf.collapse}).\qed

\begin{complement} It is worthwhile to note that condition
(\ref{top.alg.normsurfup.lem}.2) can be weakened to:
\begin{enumerate}
\item[2']  $\ex(f)$  contains a unique geometrically irreducible
surface
$S$.  $S$ has only isolated singularities and $S$ is contracted to a
point by $f$. 
\end{enumerate}
\end{complement}

 It would be very useful to have a version of 
(\ref{top.alg.normsurfup.lem}) which works if $\ex(f)$ is an
irreducible but nonnormal surface.

In the topological version (\ref{top.normsurf.collapse}) esentially
nothing can be said if $F$ is allowed to become an arbitrary compact
2-complex.  For instance, let $M$ be an arbitrary compact 3-manifold
and 
$F$  the 2-skeleton of a triangulation  of $M$. Then $\bar N$ is the
union of some copies of $S^3$ (one for each 3-simplex). 

This example usually can not arise as the real points of an algebraic
surface, but it is not hard to modify this example by approximating
each simplex with a sphere to get the following. (This is not used in
the sequel and so  no proof is given here.)

\begin{prop}\label{top.exmp}
 Let $M$ be a compact differentiable manifold of dimension
$n$. Then there is a smooth real algebraic variety $X$ and a morphism
$f:X\to Y$ with the following properties:
\begin{enumerate}
\item $X(\r)\sim M$,
\item $\overline{Y(\r)}$ is a disjoint union of copies of $S^n$,
\item $\ex(f)$ is
 a geometrically irreducible divisor and $\ex(f)(\r)$  is a union of
copies of $S^{n-1}$ intersecting transversally.
\end{enumerate}
\end{prop}

\section{The Gateway Method}

At the beginnings of the MMP, divisorial contractions were considered
to be the easily understandable part of the program and flips the hard
part. Lately, however, more and more questions  require a detailed
understanding of all the steps of the MMP.  A  fairly complete
description of all flips is known
\cite{KoMo92}, but it seems very difficult to obtain a
 list of all divisorial contractions. One can try to study the MMP in
two basic ways:

\begin{say}[Analysis of the MMP]  Starting with a  projective variety
$X$, let us run the MMP. We obtain a sequence of birational maps
$$ X=X_0\map X_1\map \cdots\map X_i\map X_{i+1}\map\cdots\map X^*.
$$ Assume that $X$ has some nice property that we would like to
preserve.  We need some way of proving  that
$X^*$ also has this property, at least under some additional
assumptions. One way is to prove this directly, by analyzing each step
of the MMP.  This would sometimes require knowing each step of the
MMP, and even in dimension 3 the list is not yet available. Still there
are many results that can be established this way, for instance the
existence of the MMP itself. In this approach one starts with a
variety $X$ and tries to understand every possible way an MMP can {\it
start} with $X$. This is oftentimes manageable if $X$ has only mild
singularities.

Another way is to look at each  step of the MMP backwards.  In
dimension 3 we have a pretty good description of the possible
singularities that arise in the course of an MMP. Thus we can start
with a variety
$Y$ and  try to understand every possible way an MMP can {\it end}
with $Y$. This also seems rather hard. Even the case when  
$Y$ is smooth is not at all understood, but in some other cases this
approach has been carried through
\cite{Kawamata9?}.  It seems that this method is  easier to apply when
$X$ is fairly singular.

The {\it gateway method} attempts to solve the original problem in an
intermediate way. In the above chain of maps there is a smallest index
$i$ such that
$X_i$ is still ``nice" but $X_{i+1}$ is not.  Hence 
$X_i\map X_{i+1}$ is a ``gateway" through which the process leaves the
set of  ``nice" varieties.  Analysing these ``gateways" should be
easier since the direct approach tends to work for the nice variety
$X_i$  and the backwards method tends to work for more complicated
singularities of
$X_{i+1}$.

Once such a list of ``gateways" is obtained,  it is a matter of
checking the list to see if some additional properties ensure that
this step does not happen.
\end{say}

One of the simplest examples where these ideas yield a nontrivial
result is the following.

\begin{exmp} Assume that we want to stay within the class of varieties
of index 1. In this case there is only one gateway:

{\bf Proposition.} {\it Let $f:X\map X'$ be a step of the
3-dimensional MMP where $X$ has index 1 but $X'$ has higher index.
Then $f$ is the contraction of a divisor
$E\subset X$ to  a point. Furthermore, 
$E\cong \p^2$,  $X$ is smooth along $E$ and $E$ has normal bundle
$\o_E(-2)$.}

This result is a special case of  \cite{Mori88} and
\cite{Cutkosky88}, though they did not approach this from the point of
view of gateways. A proof along the lines suggested by the gateway
method is not hard to construct, but this is not any shorter then the
direct proofs.

As a consequence we obtain:

{\bf Corollary} {\it  Let $X$ be a projective 3-fold with index 1
terminal singularities. Assume that $X$ does not contain any surface
$S\subset X$ which admits a birational morphism onto $\p^2$.  Then
each step of the MMP starting with $X$ is a projective 3-fold with
index 1 terminal singularities.}

Unfortunately the above condition needs to be checked for every
surface $S$, even for very singular ones. Thus in practice this does
not seem to be a useful observation.
\end{exmp}

\begin{say} Our aim is to develop a similar theory for real algebraic
threefolds. Thus we have to decide which varieties are ``nice" and
then describe all possible gateways through which the MMP can leave
the  class of ``nice" varieties.

(\ref{mmpt.concl}) naturally suggest a topological choice: $X$ is
``nice" if $X(\r)$ or $\overline{X(\r)}$ is a 3-manifold, maybe with
some additional properties. This was my first attempt, but I was
unable to make it work. The main problem seems to be that, as the
computations of \cite{rat1} show, there is basically no relationship
between the algebraic complexity of a terminal singularity $0\in X$
and the topological complexity of its real points $X(\r)$.

Eventually I settled at a completely algebraic choice: $X$ is nice if
it has index 1 along $X(\r)$. There are two  main reason for adopting
this definition:
\begin{enumerate}
\item Most complications of 3-dimensional birational geometry come
from the appearance of points of index $>1$. Hence this is likely to
be the right choice algebraically.
\item One of the first things I realized was that under this condition
there would be no flips. Indeed, flips need higher index singular
points to exist. If we have only  index 1 points along $X(\r)$, then
all higher index points appear in conjugate pairs. A look at the  
list of flips
\cite{KoMo92} shows that the singularities  appearing along a flipping
curve are {\it always} asymmetrical.
\end{enumerate}

Thus our task is to get a list of all steps $f:Y\to X$ of the MMP over
$\r$ such that $Y$ has index 1 along $Y(\r)$.  The case of
divisor--to--curve contraction is relatively easy. Most of the work is
devoted to studying the divisor--to--point contractions. Let $0\in
X(\r)$ be the  point in question.  The existence of $f$ is   local in
the Euclidean topology. I will go through the classification (up to
real analytic equivalence) of  3-dimensional terminal singularities
over $\r$  and for each describe all possible
$f:Y\to X$.

There is one subtle point here: the condition of
$\q$-factoriality is not preserved under  analytic equivalence. Thus
first we need to develop  a notion of ``extremal contraction without
$\q$-factoriality".
\end{say}

\begin{defn}\label{elem.extr.defn}
 Let $X$ be a normal variety over a field $K$ such that 
$K_X$ is $\q$-Cartier. A proper birational morphism
$f:Y\to X$ is called an {\it elementary extraction}  of
$X$ if
\begin{enumerate}
\item $Y$ is normal and $K_Y$ is $\q$-Cartier.
\item The exceptional set $\ex(f)$ contains a unique
$K$-irreducible divisor
$E$.
\item $-K_Y$ is $f$-ample.
\end{enumerate}

If we start with $Y$ and construct $f:Y\to X$ then $f$ is usually
called an {\it elementary contraction}  of $Y$.

We can write $K_Y\equiv  f^*K_X+a(E,X)E$ where $a(E,X)$ is the
discrepancy of $E$. Thus $-a(E,X)E$ is $f$-ample. An exceptional
divisor can never be relatively ample (or nef) (cf.
\cite[3.35]{KoMo98}), thus $a(E,X)>0$ and so $-E$ is
$f$-ample. This implies that 
$\ex(f)=\supp E$. 

$f(E)$ is also called the {\it center} of $f$ on $X$.
\end{defn}

A crucial property of elementary extractions is that they are
determined by their exceptional divisors:

\begin{prop}\label{gw.mats-mumf.lem} Let $X$ be a normal variety over
a field
$K$ such that 
$K_X$ is $\q$-Cartier. Let $f_i:Y_i\to X$ be  elementary extractions 
with exceptional divisors $E_i\subset Y_i$ for $i=1,2$. Assume that
$E_1$ and $E_2$ correspond to each other under the birational map
$f_2^{-1}\circ f_1:Y_1\map Y_2$.  Then $Y_1$ and $Y_2$ are isomorphic
(over
$X$).
\end{prop} 

Proof. Let $\phi:f_2^{-1}\circ f_1:Y_1\map Y_2$ be the composition.
$\phi$ is birational, and $\ex(\phi)$, $\ex(\phi^{-1})$  have
codimension at least 2. Furthermore, $K_{Y_1}$ and 
$K_{Y_2}=\phi_*(K_{Y_1})$ are relatively ample.  Thus
$\phi$ is an isomorphism by an argument of
\cite[p.671]{Matsusaka-Mumford64}.
\qed
\medskip

In some sense this gives a way of enumerating all elementary
extractions of
$X$. We try to list all exceptional divisors over $X$ and
     for each construct the corresponding unique elementary
extraction. Usually there are infinitely many elementary extractions
for a given $X$ and there does not seem to be an easy way to predict
for which divisors does the corresponding elementary extraction exist. 

The next definition singles out a special class  of
 elementary extractions, by restricting the singularities allowed on
$Y$. The aim is to formalize a special case of the gateway method: we
assume that $Y$ is ``nice".

\begin{defn}\label{gw.g-e.def}
 Let $X$ be a normal variety over a field $K$ such that 
$K_X$ is $\q$-Cartier. A proper birational morphism
$f:Y\to X$ is called a {\it gateway--extraction} or {\it
g--extraction}  if
\begin{enumerate}
\item $f$ is an elementary extraction with exceptional divisor
$E\subset Y$.
\item $Y$ has terminal singularities.
\item $K_X$ and $E$ are Cartier at the generic point of every 
geometrically irreducible $K$-subvariety of
$\ex(f)$.
\end{enumerate}
\noindent In dimension three $Y$ has only isolated singularities,
hence (3) is equivalent to the apparently weaker condition:
\begin{enumerate}
\item[3'] $K_X$ and $E$ are Cartier at every $K$-point of
$\ex(f)$.
\end{enumerate}

If we start with $Y$ and construct $f:Y\to X$ then $f$ is usually
called a {\it g--contraction}  of $Y$.
\end{defn}

The main technical aim of this article is to obtain a list of
g--extractions for threefolds with  terminal singularities.  The
project turns out to be feasible since the discrepancy 
$a(E,X)$ is always quite small. I have no a priori proof of this, but
in every case the study of low discrepancy divisors leads to a
description of all g--extractions.

The relationship between low discrepancy divisors and  g--extractions
rests on the
 following easy observation:

\begin{prop}\label{gw.d-ineq.prop}
 Let $X$ be a normal variety over a field $K$ such that 
$K_X$ is $\q$-Cartier. Let 
$f:Y\to X$ be a  g--extraction with exceptional divisor
$F\subset Y$. Let $E$ be a geometrically irreducible
$K$-divisor over $X$  such that $\cent_XE\subset
\cent_XF$. Then
$$ a(E,X)\geq a(E,Y)+a(F,X).
$$
\end{prop}

Proof. Let $g:Z\to X$ be a proper birational morphism such that
$\cent_ZE$ is a divisor on $Z$. We may assume that the induced
rational map $h:Z\map Y$ is a morphism.
$h(E)$ is a geometrically irreducible $K$-subvariety of
$Y$ which is contained in 
$\ex(f)$. Write
$$
\begin{array}{rcl} K_Z&\equiv& g^*K_X+a(E,X)E+(\mbox{other exceptional
divisors}),\\ K_Z&\equiv& h^*K_Y+a(E,Y)E+(\mbox{other exceptional
divisors}),\\ K_Y&\equiv& h^*K_X+a(F,Y)F, \qtq{and}\\ h^*F&\equiv&
cE+(\mbox{other exceptional divisors}),
\end{array}
$$ where $c>0$ since $h(E)\subset \ex(f)=\supp F$ and $c$ is an
integer by (\ref{gw.g-e.def}.3).  Making the substitutions we obtain
that
$a(E,X)= a(E,Y)+c\cdot a(F,X)\geq  a(E,Y)+a(F,X)$. \qed
\medskip

The same method also proves the following result:

\begin{prop}\label{gw.d-ineq.prop.cor}
 Let $X$ be a normal variety over a field $K$ such that 
$K_X$ is $\q$-Cartier. Let 
$f:Y\to X$ be a  morphism with exceptional divisor
$F=\cup F_i\subset Y$. Assume that $Y$ has terminal singularities and
$K_X$ and $F$ are Cartier at the generic point of every  geometrically
irreducible $K$-subvariety of
$\ex(f)$.

 Let $E$ be a geometrically irreducible
$K$-divisor over $X$  such that $\cent_XE\subset
\cup_i\cent_XF_i$. Then
$$ a(E,X)\geq a(E,Y)+\min_i\{a(F_i,X)\}.\qed
$$
\end{prop}

\begin{cor}\label{gw.discr1.cor}
 Let $X$ be a normal variety over a field $K$ such that 
$K_X$ is $\q$-Cartier.  Let 
$f:Y\to X$ be an elementary extraction with exceptional divisor
$E\subset Y$.  Assume that $E$ is geometrically irreducible and
$a(E,X)\leq 1$. 

Then either $f:Y\to X$ is a g--extraction, or $X$ has no
g--extractions whose center contains  $f(E)$.
\end{cor}

Proof.  Let $g:Z\to X$ be a g--extraction of $X$  whose center
contains  $f(E)$. Let $F\subset Z$ be the exceptional divisor. Then
$a(E,X)\geq  a(E,Z)+a(F,X)$. If
$a(E,Z)=0$ then 
$\cent_ZE$ is a divisor which is contained in $F$. Since
$F$ is an irreducible divisor, $\cent_ZE=F$, hence $Y=Z$ by
(\ref{gw.mats-mumf.lem}).  Otherwise $a(E,Y)\geq 1$ which would force
$a(F,X)\leq 0$. This contradicts (\ref{gw.g-e.def}.2).
\qed

\begin{rem} This corollary gives a very efficient way of finding all 
g--extractions of a given $X$ in some cases.  We have to find {\it
one} 
 geometrically irreducible divisor $E$ such that
$a(E,X)\leq 1$ and construct the corresponding elementary extraction
$f:Y\to X$. Then it is usually easy to determine the singularities of
$Y$.

\cite{Markushevich96} proved that if $0\in X$ is a terminal threefold
singularity which is not smooth, then there is a divisor $E$  over
$\bar K$ with
$\cent_XE=\{0\}$ and $a(E,X)\leq 1$. Thus there is always such an
irreducible $K$-divisior, but it may not be geometrically irreducible.
Still, in many cases we are  able to apply (\ref{gw.discr1.cor})
directly.

In the remaining cases we show that there is always a  geometrically
irreducible divisor $E$   with
$\cent_XE=\{0\}$ and $a(E,X)\leq 3$. This is still very useful, thanks
to the following:
\end{rem}

\begin{cor}\label{gw.discr2.cor}
 Let $X$ be a normal variety over a field $K$ such that 
$K_X$ is Cartier.  Let 
$f:Y\to X$ be an elementary extraction with exceptional divisor
$E\subset Y$.  Assume that $E$ is geometrically irreducible.  Let
$g:Z\to X$ be any g--extraction with exceptional divisor $F$ whose
center contains  $f(E)$.

Then either $g=f$ or $a(F,X)\leq a(E,X)-1$. 
\end{cor}

Proof.  By (\ref{gw.d-ineq.prop}),  $a(E,X)\geq  a(E,Z)+a(F,X)$. If
$a(E,Z)=0$ then $E$ and
$F$ correspond to each other, hence $Y=Z$ by
(\ref{gw.mats-mumf.lem}).  Otherwise $a(E,Y)\geq 1$, thus
$a(F,X)\leq a(E,X)-1$.  \qed

\section{Small and Divisor--to--Curve Contractions}

In this section we look at those steps $f:X\to Y$ of the MMP over $\r$
which are either small contractions or contract a divisor to a curve.
The two cases  can be treated together  in the following setting:

\begin{notation}\label{sm.not}
 Let $K$ be a field of characteristic 0. Let $X$ be a 3-fold over $K$
with terminal singularities  and $f:X\to Y$  a proper birational
morphism over $K$ such that $-K_X$ is $f$-ample and $f_*\o_X=\o_Y$.
Let $0\in Y(K)$ be a closed point such that $\dim f^{-1}(0)=1$. 

We will need that under these assumptions 
$R^1f_*\o_X=R^1f_*\o_X(K_X)=0$ by the generalized
Grauert--Riemenschneider vanishing theorem (see, for instance, 
\cite[8.8]{CKM88} or \cite[2.65]{KoMo98}). 
\end{notation} 

In keeping with the principles of the gateway method, we are interested
in the case when
$X$ has index 1 at all points of $X(K)$. The following theorem gives a
complete description of such contractions:

\begin{thm}\label{sm.thm}
 Notation and assumptions as in (\ref{sm.not}).  Assume in addition 
that $X$ has index 1 at all points of $X(K)$. Then $Y$ is smooth at
$0$ and one can choose local  (analytic or formal) coordinates 
$(x,y,z)$ at $0\in Y$ such that $X$ is the blow up of the curve
$(z=g(x,y)=0)\subset Y$ for some $g\in K[[x,y]]$.

In particular, $f$ can  not be small. 
\end{thm} 

This theorem has some very useful consequences for the MMP over $\r$:

\begin{cor}\label{sm.mmp.cor1}   Starting with a  projective variety 
$X$ over $\r$, let 
$$ X=X_0\map X_1\map \cdots\map X_{i}\stackrel{f_{i}}{\map} X_{i+1}
$$  be the beginning of an  MMP over $\r$.  Assume that   $X_j$ has
index 1 at all points of $X_j(\r)$ for
$j\leq i$.  Then the induced  maps between the sets of real points
$$ X(\r)=X_0(\r)\to  \cdots\to X_{i}(\r)\stackrel{f_{i}}{\to}
X_{i+1}(\r)
$$  are  everywhere defined.
\end{cor}

Proof. The only steps of the MMP over $\r$ which are not everywhere
defined are the flips of small contractions (\ref{mmpt.cases.prop}). 
By (\ref{sm.thm}) there are no flips in the sequence.\qed
\medskip

The topological behavior of divisor--to--curve contractions can also
be determined using (\ref{sm.thm}):

\begin{thm}\label{sm.top.thm} Let $X$ be a proper 3-fold over $\r$
with terminal singularities  such that $X$ has index 1 at all points of
$X(\r)$ and $\overline{X(\r)}$ is a 3-manifold. Let $f:X\to Y$ be  a
proper birational morphism over $\r$ such that $-K_X$ is $f$-ample and
$f_*\o_X=\o_Y$. Assume that $\dim f^{-1}(y)\leq 1$ for every $y\in
Y$.   Then either
\begin{enumerate}
\item $f$ is $\r$-small, or
\item $\overline{X(\r)}$ contains a 1-sided torus or Klein bottle with
nonorientable neighborhood. 
\end{enumerate}
\end{thm}

Proof.  By (\ref{sm.thm}), there is a real curve $D\subset Y$ such
that 
$Y$ is smooth along $D$ and $X=B_DY$ (at least in a neighborhood of
$Y(\r)$).   Pick $0\in D(\r)$  and let $(z=g(x,y)=0)$ be a local
equation of
$D$. By  (\ref{sm.ci.blowup}),  either $D$ is smooth at $0$ or 
$X$ has a unique singular point over $0$ with local equation
$st=g(x,y)$, which is equivalent to
$s^2-t^2-g(x,y)=0$. These are of type
$cA_{>1}^-$ or $cA_1$ in the classification of \cite{rat1}. If $g$
does not change sign on the $(x,y)$-plane then 
$X(\r)\setminus f^{-1}(0)\to Y\setminus\{0\}$ is one--to--one near $0$,
hence
$f$ is $\r$-small near $0$. If $g$  does change sign on the
$(x,y)$-plane, then from \cite[sec. 4]{rat1} we see that (after a
coordinate change) $g=\pm(x^2+y^{2r+1})$ and
$X(\r)$ is a manifold near $f^{-1}(0)$.  In particular, $D(\r)$ is the
disjoint union of some isolated points and some copies of $S^1$. 

If $D(\r)$ is finite then
$f$ is $\r$-small. Otherwise $D(\r)$ has a connected component
 $M\sim S^1$. Let $E\subset X$ be the exceptional divisor of
$f$. By explicit computation we see that 
$E(\r)\to D(\r)$ is an $S^1$-bundle. Hence there is a unique connected
component $N\subset E(\r)$ such  that  $N$ is an $S^1$-bundle over
$M$. Thus
$N$ is either a torus or a Klein bottle. $N$ is 1-sided with
nonorientable neighborhood, since these hold locally for the blow up of
a smooth curve in a smooth 3-fold.\qed

\begin{exmp}\label{sm.ci.blowup}
 Set $Y={\Bbb A}^3$ with
 coordinates 
$(x,y,z)$. Let $X$ be the blow up of the curve
$(z=g(x,y)=0)\subset Y$. Then $X$ has a unique singular point which is
given by an equation $st-g(x,y)=0$. 
\end{exmp}

\begin{cor}\label{sm.mmp.cor2}
  Starting with a  projective variety 
$X$ over $\r$, let 
$$ X=X_0\map X_1\map \cdots\map X_{i}\stackrel{f_{i}}{\map} X_{i+1}
$$  be the beginning of an  MMP over $\r$.  Assume that 
\begin{enumerate}
\item  $X_j$ has index 1 at all points of $X_j(\r)$ for
$j\leq i$,  
\item  $\overline{X_j(\r)}$  is a PL-manifold for
$j\leq i$,  
\item $\overline{X(\r)}$ satisfies the conditions (\ref{int.no.cond}).
\end{enumerate}
\noindent Then: 
\begin{enumerate}
\setcounter{enumi}{3}
\item The induced  maps between the sets of real points
$f_j: X_{j}(\r)\to X_{j+1}(\r)$  are  everywhere defined   for  $j\leq
i$,
\item  For every $j\leq i+1$, there is a finite set
$S_j\subset X_{j}(\r)$ such that
$\overline{X_{j}(\r)}\setminus S_j$ is homeomorphic to  an open subset
of  $\overline{X(\r)}$, 
\item The smooth part of $\overline{X_{i+1}(\r)}$  also satisfies the
conditions (\ref{int.no.cond}).
\end{enumerate}
\end{cor}

Proof. The  steps of an MMP  are  everywhere defined  by
(\ref{sm.mmp.cor1}). 

If $g:U\to V$ is any divisorial contraction over $\r$ then
$\ex(g^{-1})(\r)$ is finite unless $g$ is a  divisor--to--curve
contraction  which is not $\r$-small. 

Let $f_j$ be the first divisor--to--curve contraction in the sequence
which is not $\r$-small. By the above remark, (5) holds for
$j$. By (\ref{sm.top.thm}),
$X_j(\r)$ contains a  surface $F$ which is either a 1-sided torus or
Klein bottle with nonorientable neighborhood. We can move $F$ away from
any finitely many points, thus by (5) 
$X(\r)$ also contains a   1-sided torus or Klein bottle with
nonorientable neighborhood. This is a contradiction. Hence among the
steps  there is no  divisor--to--curve contraction  which is not
$\r$-small. This gives (5) and (6).\qed

\medskip

The proof of (\ref{sm.thm}) relies on  two results:

\begin{prop}\cite[Thm. 4]{Cutkosky88} \label{sm.cut}  (\ref{sm.thm})
holds if $K$ is algebraically closed and $C$ is irreducible.\qed
\end{prop}

\begin{lem}\label{sm.tree}(cf. \cite[1.14]{Mori88})    Let $f:X\to Y$
be a proper morphism and 
$0\in Y$  a closed point such that $\dim f^{-1}(0)=1$. Set $\red
f^{-1}(0)=C=\cup C_i$. 
\begin{enumerate}
\item If $R^1f_*\o_X=0$ then $C$ is a tree of smooth rational curves.
\item Let $D$ be a $\q$-Cartier Weil divisor on $X$ such that $D$ is
Cartier at all but finitely many points of $f^{-1}(0)$. Assume that
$(D\cdot C_i)<0$ for every $i$ and $R^1f_*\o_X(D)=0$. Then 
$-1\leq (D\cdot C_i)<0$ for every $i$ and
$D$ is not Cartier at the singular points of $C$.
\end{enumerate}
\end{lem}

Proof. By replacing $Y$ with a neighborhood of $0$, we may assume that
every fiber of $f$ has dimension at most 1. 

Let $G$ be a sheaf on $X$ such that $R^1f_*G=0$ and $Q=G/F$ a quotient
of
$G$ whose support is in $f^{-1}(0)$. We get an exact sequence
$$ R^1f_*G\to R^1f_*Q\to R^2f_*F.
$$ The left hand side is zero by assumption and the right hand side is
zero since every fiber of $f$ has dimension at most 1. Thus
$R^1f_*Q=0$.

Applying this with $G=\o_X$ and $Q=\o_C$ we conclude that 
$H^1(C,\o_C)=R^1f_*\o_C=0$, hence $C$ is a tree of smooth rational
curves. This proves (1).

In order to see the second part, we may assume that the residue field
of
$0$ is algebraically closed. Then a point   $P\in C$ is singular iff
there are at least 2 irreducible components through  $P$. 

$\o_X(D)\otimes \o_{C_i}$ is a rank one locally free sheaf except
possibly at the ponts where $D$ is not Cartier. Let $L_i$ denote its
quotient by the torsion subsheaf. Then $L_i$ is an invertible sheaf
and we have a surjection
$\o_X(D)\to L_i$. Applying $R^1f_*$ we obtain as above that
$H^1(C_i,L_i)=0$. Thus $\deg L_i\geq -1$.

On the other hand, for every
$m>0$ we have an injection
$$ L_i^m\cong  (\o_X(D)^{\otimes m}\otimes \o_{C_i})/(\mbox{torsion})
\into  (\o_X(mD)\otimes \o_{C_i})/(\mbox{torsion}).
$$ If $mD$ is Cartier  then the right hand side has negative degree,
thus
$L_i^m$ has negative degree. Therefore $\deg L_i=-1$ for every $i$.
Furthermore, $m(D\cdot C_i)\geq m\deg L_i=-m$, so 
$(D\cdot C_i)\geq -1$.

Set $M:=(\o_X(D)\otimes \o_{C})/)\mbox{torsion})$. 
$H^1(C,M)=0$ as above. We have an exact sequence
$$ 0\to M\to \sum L_i\to Q\to 0,
$$ where $Q$ is supported at the singular points of $C$. Taking
cohomologies, we conclude that $H^0(C,Q)=0$ thus $Q=0$. 

If $D$ is Cartier at a singular point $P$ of $C$ then $M$ is locally
free at $P$ and $M\to \sum L_i$ can not be surjective at
$P$ (it is not even surjective when tensored with the residue field at
$P$).\qed

\begin{cor}\label{sm.tree2}
 Notation and assumptions as in (\ref{sm.not}). Then
 $C$ is a tree of smooth rational curves and  $K_X$ is not Cartier at
the singular points of $C$.
\end{cor}

Proof. Apply (\ref{sm.tree}) with $D=K_X$.
$R^1f_*\o_X=R^1f_*\o_X(K_X)=0$ by (\ref{sm.not}). \qed

\begin{say}[Proof of (\ref{sm.thm})]
 The assumptions and conclusions are local near
$0$, thus we may replace $Y$ by a suitable analytic or formal
neighborhood of
$0$.

By (\ref{sm.tree2}),  $C$ is a connected tree of smooth rational
curves. $\gal(\bar K/K)$ acts on $C$, thus $C$ either has  a
 singular $K$-point or a geometrically irreducible component defined
over
$K$.

If $P\in C$ is a singular point then $K_X$ is not Cartier at $P$ by
(\ref{sm.tree2}), but if $P\in X(K)$ then $K_X$ is Cartier at $P$ by
assumption. Thus $C$ can not have a  singular $K$-point.

Let $C_0\subset C$ be a geometrically  irreducible component defined
over
$K$. Let $H\subset X$ be a divisor defined over $K$ which intersects
all irreducible components of $C\setminus C_0$ transversally but is
disjoint from $C_0$.  A large multiple of $H$ defines a morphism $X\to
Y'\to Y$ such that
$C_0$ is contracted to a point in $Y'$.  If (\ref{sm.thm}) holds for
$X\to Y'$, then $X$ has index one along $C_0$. By (\ref{sm.tree2}) this
implies that $C_0$ is a connected component of 
$C$. On the other hand, $C$ is connected since $f_*\o_X=\o_Y$. Thus
$C=C_0$ and $Y'=Y$.

Therefore it is sufficient to prove (\ref{sm.thm}) under the additional
assumption that $C$ is geometrically irreducible.

First we show that (\ref{sm.thm}) holds if $X$ has only index 1 points
along
$C$. By (\ref{sm.cut}), $Y$ is smooth at $0$ and $X=B_DY$ where
$D\subset Y$ is a  curve of embedding dimension 2.
$D$ is the image of the exceptional divisor of $f$, hence $D$ is
defined over $K$. Since 
$D$ has embedding dimension 2, its ideal is of the form $(z,g(x,y))$.

Finally we show that $X$ has only index 1 points along $C$.  We start
with the case when $K=\r$. Let $P_1,\bar P_1, \dots, P_k,\bar P_k$ be
all the conjugate pairs of points of index $>1$.  At each $P_i$ pick a
local member  $D_i\in |K_X|$  such that $C\cap D_i=P_i$. (In order to
do this, we may need to replace $X_{\bar K}$ with a smaller analytic 
neighborhood of $C$.) Let $\bar D_i$ be the conjugates. Set
$D=\sum D_i$. Let $m>1$ be the smallest natural number such that $mD$
is Cartier.
$D-\bar D$ is a Weil divisor and $\o_X(m(D-\bar D))\cong \o_X$  since
the Picard group of a neighborhood of $C$ is isomorphic to
$H^2(C(\c),\z)$ (cf. \cite[4.13]{KoMo98}). 

Corresponding to $1\in H^0(X,\o_X)$ we obtain an $m$-sheeted cyclic
cover $\pi: \tilde X\to X$ which is unramified outside the points of
index
$>1$.  Thus $K_{\tilde X}=\pi^*K_X$ and $\tilde X$ has index 1 terminal
singularities. Let $\tilde f:\tilde X\to \tilde Y$ be the Stein
factorization of $\tilde X\to Y$. By the already discussed index 1
case,
$\tilde Y$ is smooth and one can choose local analytic coordinates 
$(x,y,z)$ at $0\in \tilde Y$ such that $\tilde X$ is the blow up of the
curve
$(z=g(x,y)=0)\subset \tilde Y$. 

The group  of $m^{th}$ roots of unity (denoted by $\z_m$)  acts on
$\tilde f:\tilde X\to \tilde Y$ and the quotient is $ f: X\to  Y$. If
$\mult_0g\geq 2$ then $\tilde X$ has a unique singular point
(\ref{sm.ci.blowup}), which is necessarily fixed by the
$\z_m$-action.  Thus $X$ would have a unique point (of index $m$)
which is the quotient of a singular point.  On the other hand, the
index $>1$ singularities of
$X$ come in conjugate pairs. Therefore
$\tilde X$ is smooth and $\tilde f:\tilde X\to \tilde Y$ is the blow up
of a smooth curve $(z=y=0)\subset \tilde Y$. 

We can choose local coordinates $(x,y,z)$ on $\tilde Y$ such that the
action is
$$ (x,y,z)\mapsto (\epsilon^a x,\epsilon^by, \epsilon^cz)
$$ where $\epsilon$ is a primitive $m^{th}$ root of unity.  The
corresponding action on  $\tilde X$  has two fixed points (or a fixed
curve)  and the corresponding quotients are
$$
\c^3/{\textstyle \frac1{m}}(a,b-c,c)\qtq{and}
\c^3/{\textstyle \frac1{m}}(a,b,c-b).
$$ These are both of type $cA_0/n$ on the list (\ref{mmp.ts.thm}). A
simple checking  shows that both of these can not be simultaneously
terminal.

If $K$ is arbitrary, we can still proceed as above if we can find local
divisors $D_i\in |K_X|$ at the index $>1$ points such that $(C\cdot
\sum D_i)=0$. Finding the $D_i$ needs a little case by case analysis,
and sometimes it can be done only after first taking an auxiliary
cover. It is probably easier to observe that there can be at most 2
points of index $>1$ along $C$  (see, for instance,
\cite[14.5.5]{CKM88}), thus in  fact the only case we need to handle
is when there is precisely one pair of conjugate points of index $>1$.
\qed
\end{say}

\section{Proof of the Main Theorems}

The determination of all divisor--to--point g--extractions is rather
technical and lengthy. In this section I state
 a summary of the list of all g--extractions, and then use it to prove
the main theorems stated in the introduction.  The proofs of 
(\ref{ge.gwsing.thm}) and of (\ref{ge.gwextr.thm}) are given in
sections 9--11.

\begin{notation}\label{7.notation}
 Let $g(x_1,\dots,x_m)$ be a polynomial or power series and let $M$ be
a monomial in the $x_i$.   $M\in g$
 means that $M$ appears in $g$ with nonzero coefficient.
\end{notation}

\begin{thm}\label{ge.gwsing.thm}
 Let $0\in X$ be a three dimensional terminal singularity over
$\r$. If $X$ has a g--extraction then $0\in X$ is one of the
  following (up to real analytic equivalence near $0$).
\begin{enumerate}

\item ($cA_0$) Smooth point.

\item ($cA_0/2$) Quotient of a smooth point by the
$\z_2$-action 
$(x,y,z)\mapsto (-x,-y,-z)$.

\item ($cA_{>0}^+$) Given as
 $(x^2+y^2+g_{\geq m}(z,t)=0)$, where $g_m(z,t)\neq 0$ and
$m\geq 2$. 

\item ($cA_{>0}^+/2$) Given as
 $(x^2+y^2+g_{\geq m}(z,t)=0)/\z_2$, where $g_m(z,t)\neq 0$, $m\geq 2$
and the
$\z_2$-action is $(x,y,z,t)\mapsto (-x,-y,-z,t)$.
 Furthermore, one of the following two conditions has to be satisfied:
\begin{enumerate}
\item $m$ is divisible by $4$ and $z^m,t^m\in g$, or
\item\label{7.bad.index.2} $m$ is odd. 
\end{enumerate}
\item ($cE_6$)  Given as
$(x^2+y^3+(z^2+t^2)^2+ yg_{\geq 4}(z,t)+h_{\geq 6}(z,t)=0)$.
\end{enumerate}
\end{thm}

\begin{complement} All the singularities on the above list have
g--extractions, with the possible exception  of types
(\ref{7.bad.index.2}).  These singularities have not been analyzed
completely.
\end{complement}

\begin{thm}\label{ge.gwextr.thm}
 Let $0\in X$ be a three dimensional terminal singularity over
$\r$ and $f:Y\to X$   a g--extraction with exceptional divisor
$E=\red f^{-1}(0)$. If  
$E$ is geometrically irreducible then $f:Y\to X$ is on the following
list (up to real analytic equivalence near $0$).
\begin{enumerate}

\item ($cA_0$, point blow up)  $B_0{\Bbb A}^3\to {\Bbb A}^3$,
$E\cong
\p^2$.

\item ($cA_0$, curve blow up)  $B_C{\Bbb A}^3\to {\Bbb A}^3$ where
$C\subset {\Bbb A}^3$ is a geometrically irreducible, real and
locally  planar curve.

\item ($cA_0/2$) $B_0{\Bbb A}^3/\z_2\to {\Bbb A}^3/\z_2$, where the
$\z_2$-action on ${\Bbb A}^3$ is
$(x,y,z)\mapsto (-x,-y,-z)$. 
$E\cong \p^2$.

\noindent Furthermore, in this case there are no other g--extractions
whose center contains the origin.

\item ($cA_{>0}^+$, $\mult_0g$ even) 
 $X=(x^2+y^2+g_{\geq 2m}(z,t)=0)$ where $g_{2m}(z,t)\neq 0$ and $m\geq
1$. 
$Y=B_{(m,m,1,1)}X$
 and 

$E=(x^2+y^2+g_{2m}(z,t)=0)\subset \p^3(m,m,1,1)$. 

\noindent Furthermore, in this case there are no other g--extractions
whose center  contains  the origin.

\item ($cA_{>0}^+$, $\mult_0g$ odd) 
 $X=(x^2+y^2+g_{\geq 2m+1}(z,t)=0)$ where $g_{2m+1}(z,t)\neq 0$ and
$m\geq 1$. This case occurs only if  there is a linear change of the
$(z,t)$-coordinates  such that 
$t^{2m+1}\in g$ and 
$z^it^j\not\in g$ for  $i+2j< 4m+2$. In this coordinate system,    
$Y=B_{(2m+1,2m+1,1,2)}X$ and

$E=(x^2+y^2+g_{2m+1}(z,t)=0)\subset \p^3(2m+1,2m+1,1,2)$.

\item ($cA_{>0}^+/2$, $\mult_0g$ even)  
 $X=(x^2+y^2+g_{\geq 2m}(z,t)=0)/\z_2$, where $g_{2m}(z,t)\neq 0$,
$m\geq 1$ and the
$\z_2$-action is $(x,y,z,t)\mapsto (-x,-y,-z,t)$. This case occurs
only if $m$ is even and $z^{2m},t^{2m}\in g$. Then
$Y=B_{(m,m,1,1)}\tilde X/\z_2$ and $E=\tilde E/\z_2$, where 

\noindent
$\tilde X=(x^2+y^2+g_{\geq 2m}(z,t)=0)$ and

\noindent
$\tilde E=(x^2+y^2+g_{2m}(z,t)=0)\subset \p^3(m,m,1,1)$. 

Furthermore, in this case there are no other g--extractions whose
center contains  the origin.

\item ($cA_{>0}^+/2$, $\mult_0g$ odd) 
 $X=(x^2+y^2+g_{\geq 2m+1}(z,t)=0)/\z_2$ where $g_{2m+1}(z,t)\neq 0$,
$m\geq 1$  and the
$\z_2$-action is $(x,y,z,t)\mapsto (-x,-y,-z,t)$. In this case  I do
not have a complete list.
\end{enumerate}
\end{thm}

\begin{cor}\label{ge.Enorm.cor}
 Let $0\in X$ be a three dimensional terminal singularity over
$\r$ and $f:Y\to X$   a g--extraction with exceptional divisor
$E=\red f^{-1}(0)$.  Assume that we are not in case
(\ref{ge.gwextr.thm}.7).  If $E$ is geometrically irreducible then
$E$ is normal.
\end{cor}

Proof.  Equations for $E$ are given in (\ref{ge.gwextr.thm}).
$E\cong
\p^2$ in the first two cases. In the remaining cases  $E$ is (or is
the quotient of) a surface of the form
$$ F:=(x^2+y^2+p(z,t)=0)\subset \p^3(r,r,1,s).
$$
 All the singularities of $F$ are contained  in the
$(x=y=0)$ line. Thus we get only finitely many singularities if $p$ is
not identically zero, which is always the case in
(\ref{ge.gwextr.thm}).\qed

\begin{say}[Proof of (\ref{ge.gwsing.thm}) and (\ref{ge.gwextr.thm})
$\Rightarrow$  (\ref{int.mmp.sings}) and (\ref{int.mmp.steps})]{\ }

Under the additional assumption that $X(\r)$ satisfies the conditions
(\ref{int.no.cond}), we need to exclude the cases
(\ref{ge.gwsing.thm}.2),  (\ref{ge.gwsing.thm}.4),
(\ref{ge.gwextr.thm}.2) and in (\ref{ge.gwsing.thm}.3) we need to show
that  $g$ is not everywhere negative in a punctured neighborhood of
$0$.

Starting with $X$, let us run the MMP over $\r$. We get a sequence
$$ X=X_0\map X_1\map \cdots\map X_{i}\stackrel{f_{i}}{\map} X_{i+1}.
$$  Assume by induction that  (\ref{int.mmp.sings})  holds for $X_j$
for $j\leq i$ and (\ref{int.mmp.steps}) holds for 
$f_j:X_j\map X_{j+1}$ for $j\leq i-1$. 

We need to show that (\ref{int.mmp.sings})  holds for $X_{i+1}$ and
(\ref{int.mmp.steps}) holds for 
$f_i:X_i\map X_{i+1}$.

By (\ref{sm.mmp.cor2}) $f_j:X_j(\r)\to X_{j+1}(\r)$ are  everywhere
defined  for $j\leq i-1$ and $\overline{X_{i}(\r)}$  does not contain
a 1-sided torus or Klein bottle with nonorientable neighborhood.
Furthermore, by (\ref{sm.top.thm}), $f_i$ is  either  $\r$-small or a
divisor--to--point contraction.

By induction
$X_i$ has index 1 along $X_i(\r)$, thus $f_i$ is a g--extraction,
hence it is one of the cases listed in (\ref{ge.gwsing.thm}). We
excluded several cases one at  a time.

\medskip {\it Excluding (\ref{ge.gwsing.thm}.2).} By
(\ref{ge.gwextr.thm}.3)
$X_i\to X_{i+1}$ is the blow up of the singular point $\a^3/\z_2$.
This gives a 1--sided $\r\p^2$ in $X_i(\r)$. This is  a contradiction
by (\ref{sm.mmp.cor2}).

\medskip {\it Excluding (\ref{ge.gwextr.thm}.2).} In this case
$X_i$ contains a 1-sided torus or Klein bottle with nonorientable
neighborhood by (\ref{sm.top.thm}). This is again a contradiction by
(\ref{sm.mmp.cor2}).

\medskip {\it Excluding (\ref{ge.gwsing.thm}.3) with $g<0$.}  That is,
we  consider the case when $0\in X_{i+1}$ is of the form
$(x^2+y^2+g_{\geq m}(z,t)=0)$ and  $g$ is  everywhere negative in a
punctured neighborhood of $0$. (These are called $cA^+_{>0}(0,-)$ in
\cite{rat1}.) By \cite[4.4]{rat1} the link of $0\in X_{i+1}(\r)$ is a
torus. This gives only a 2-sided torus in $X(\r)$ which is allowed. I
proceed to prove, however, that    we still get a 1-sided torus in
$X(\r)$ coming from   the exceptional divisor of
$X(\r)\to X_{i+1}(\r)$. This   contradicts (\ref{int.no.cond}.2).

 $m=\mult_0g(z,t)$  is necessarily even, say
$m=2r$. By (\ref{ge.gwextr.thm}.4) the only g--extraction is the
$(r,r,1,1)$-blow up.  Thus $X_i\to X_{i+1}$ is this blow up. We
distinguish two cases:

General case: $g_{2r}(z,t)$ is
  negative on $\r^2\setminus\{0\}$. The exceptional divisor $E$ of the
above g--extraction  is the weighted hypersurface
$$ E=(x^2+ y^2+g_{2r}(z,t)=0)\subset \p(r,r,1,1).
$$ Its canonical divisor is $K_E=\o_E(-2)$, thus $E$ is orientable.
The projection $(x:y:z:t)\mapsto (z:t)$ exhibits $E$ as an
$S^1$-bundle over $\r\p^1$, thus
$E\sim S^1\times S^1$. $L(0\in X_{i+1}(\r))$ is connected, thus
$E(\r)\subset X_i(\r)$ is a 1-sided torus. By (\ref{sm.mmp.cor2}), 
$X(\r)$ also contains a 1-sided torus, a contradiction.

Special case: $g_{2r}(z,t)$ is not
 negative on $\r^2\setminus\{0\}$. 
$g_{2r}(z,t)$ is the leading term of 
$g_{\geq m}(z,t)$, which  is 
 negative on $\r^2\setminus\{0\}$. Thus 
$g_{2r}(z,t)$ is nonpositive
  on $\r^2\setminus\{0\}$. 

The $t$-chart on $X_i\cong B_{(r,r,1,1)}X_{i+1}$ is
$x_1^2+ y_1^2+ t_1^{-2m}g(z_1t_1,t_1)$. Set
$g'(z_1,t_1):=t_1^{-2m}g(z_1t_1,t_1)$. Then 
$g'(z_1,t_1)$ is strictly negative outside the
$z_1$-axis, and is not identically zero on the
$z_1$-axis. Thus $g'(z_1,t_1)$ is everywhere nonpositive with only
finitely many zeros. Thus at each zero of
$g'(z_1,t_1)$, $X_i$ has a singular point of type
$cA^+_{>0}(0,-)$. This contradicts the inductive assumption.

Thus $X_{i+1}$ does not contain any  
 $cA^+_{>0}(0,-)$ type points.

\medskip {\it Excluding (\ref{ge.gwsing.thm}.4).} 
\cite[5.9]{rat1} shows that  the link of a singularity of type
$cA^+/2$ contains a connected component homeomorphic to
$\r\p^2$, except when we can write the singularity as
$(x^2+y^2+g_{\geq m}(z,t)=0)/\z_2$, where the
$\z_2$-action is $(x,y,z,t)\mapsto (-x,-y,-z,t)$ and $g$ is 
everywhere negative in a punctured neighborhood of $0$.

If the link of a singularity of $X_{i+1}(\r)$ 
 contains a connected component homeomorphic to
$\r\p^2$, then  $X_i(\r)$ contains a 2-sided 
$\r\p^2$. Hence by (\ref{sm.mmp.cor2}), 
$X(\r)$ also contains a 2-sided 
$\r\p^2$ which is excluded by (\ref{int.no.cond}.1).

Thus we are reduced to the case
$(x^2+y^2+g_{\geq m}(z,t)=0)/\z_2$, where the
$\z_2$-action is $(x,y,z,t)\mapsto (-x,-y,-z,t)$ and $g$ is 
everywhere negative in a punctured neighborhood of $0$. As in the
previous case,  $m=\mult_0g(z,t)$  is necessarily even, say
$m=2r$ and by (\ref{ge.gwextr.thm}.6) the only g--extraction is the
$(r,r,1,1)$-blow up. We proceed to  prove that   the exceptional
divisor of
$X_i(\r)\to X_{i+1}(\r)$ contains a 1-sided Klein bottle with
nonorientable neighborhood, contradicting (\ref{int.no.cond}.3). We
again distinguish two cases:  

General case: $g_{2r}(z,t)$ is
 negative on $\r^2\setminus\{0\}$. The exceptional divisor $E$ is the
$\z_2$ quotient of the weighted hypersurface
$$
\tilde E=(x^2+ y^2+g_{2r}(z,t)=0)\subset \p(r,r,1,1).
$$  We already  determined that  $\tilde E(\r)$ is a torus and
 a choice of orientation is given by
$(dy\wedge dz)/x$. The $\z_2$-action sends this to
$(d(-y)\wedge d(-z))/(-x)=-(dy\wedge dz)/x$, hence
$\tilde E(\r)/\z_2$ is not orientable. We conclude that one of the
connected components of $E(\r)$ is a Klein bottle.  (There may be
other connected components.) The Klein bottle is 1-sided since 
$\tilde E(\r)$ is 1-sided. The regular 
 neighborhood is nonorientable  since its boundary, the link of $0\in
X(\r)$, is again a Klein bottle.

Special case: $g_{2r}(z,t)$ is not negative on $\r^2\setminus\{0\}$. 

The same computation as above shows that  this leads to a $cA^+/2$ 
point of the same type that we started with on
$X_i$, which contradicts the inductive assumption.

Thus we conclude that $X_{i+1}$ does not contain any  
$cA^+/2$   type points.\qed
\end{say}

\begin{say}[Proof of (\ref{int.orient.thm}) and
(\ref{int.nonorient.thm})]{\ }

We follow the steps of an MMP over $\r$, using  (\ref{int.mmp.steps}). 
$f_i:X_i(\r)\to X_{i+1}(\r)$ is a homeomorphism in cases
(\ref{int.mmp.steps}.1--2) while (\ref{int.mmp.steps}.3) gives a
connected sum with $\r\p^3$.

In the cases (\ref{int.mmp.steps}.4) the exceptional divisor is normal
by (\ref{ge.Enorm.cor}), hence we get various cases of
(\ref{int.orient.thm}) by (\ref{top.alg.normsurfup.lem}). \qed
\end{say}

\begin{exmp}\label{ge.2extr.exmp}  Consider the singularity
$X:=(x^2+y^2+z^{2m+1}+t^{4m+2}=0)$. The 
$(2m+1,2m+1,2,1)$ blow up   $X_1\to X$ is  a g--extraction which is
smooth along the $\r$-points.

The $(m,m,1,1)$-blow up is  another g--extraction whith one singular
point
$(x_1^2+y_1^2+z_1^{2m+1}t_1+t_1^{2m+2}=0)$ on the
$t$-chart. After the $(m+1,m+1,1,1)$-blow up we obtain a variety
$X_2\to X$ which is smooth along its $\r$-points.

These two resolutions are indeed quite different. Using the methods of
section 5, we see  that
$X_1(\r)\sim X(\r)\ \#\ \r\p^3$ and
$X_2(\r)\sim X(\r)\ \#\ S^1\times S^2$.
\end{exmp}

\section{$cAx$ and $cD$-type Points}

In this section we begin to classify g--extractions (\ref{gw.g-e.def})
of terminal singularities over any field.  The classification of
3--fold terminal singularities over nonclosed fields is done in
\cite{rat1}. The results are summarized in (\ref{mmp.ts.thm}). We work
through the list of the singularities. In most cases it is easy to see
that there are no g--extractions. This is done by exhibiting an
elementary extraction (\ref{elem.extr.defn}) which is not a
g--extraction. If the discrepancy   of the exceptional divisor is
$\leq 1$ then there are no g--extractions by (\ref{gw.discr1.cor}).

In this section we deal with the cases $cAx/2, cAx/4, cD, cD/2, cD/3$.
Among terminal singularities these are somewhat esoteric but the
proofs work well for them: in each case (\ref{gw.discr1.cor}) applies.

The remaining terminal singularities are considered in the next 2
sections. In some cases much more complicated arguments are needed to
classify all g--extractions.

\begin{defn}\label{8.wbup}[Weighted blow-ups]{\ }

Let $x_1,\dots,x_n$ be coordinates on $\a^n$. The  usual blow up of
the origin is patched together from affine charts with morphisms of
the form
$$ x_j=x'_jx'_i\qtq{if} j\neq i\qtq{and}  x_i=x'_i.
$$ I   refer to this as the {\it $x_i$-chart}.

Let $a_1,\dots,a_n$ be  a sequence of positive integers. For every
$1\leq i\leq n$  we can define a morphism  $\Pi_i:\a^n\to \a^n$ by
$$ x_j=x'_j(x'_i)^{a_j}\qtq{if} j\neq i\qtq{and}  x_i=(x'_i)^{a_i}.
$$ This morphism  is birational iff $a_i=1$ and has degree $a_i$ in
general. One can easily notice that $\Pi_i$ is invariant under the
action
$$
\a^n(x'_1,\dots,x'_n)/\textstyle{\frac1{a_i}}
(-a_1,\dots,-a_{i-1},1,-a_{i+1},\dots,-a_n)
$$ and it descend to  a birational morphism $\pi_i$
$$
\Pi_i:\a^n(x'_1,\dots,x'_n)\to \a^n(x'_1,\dots,x'_n)/\z_{a_i}
\stackrel{\pi_i}{\longrightarrow} \a^n(x_1,\dots,x_n).
$$ Furthermore, these charts patch together to give a projective
morphism 
$$
\pi:B_{(a_1,\dots,a_n)}\a^n\to \a^n.
$$ This is called the  {\it weighted blow up} of $\a^n$ with weights
$a_1,\dots,a_n$.
\end{defn}

\begin{notation}\label{8.notation}
 In the proofs in sections 9--11 I use the following conventions.

Firts I state the name of the singularity $X$ from (\ref{mmp.ts.thm})
and possibly some other  restrictions. Then I write down the normal
form of the equation   $X=F(x,y,z,t)/\frac1{r}(b_x,b_y,b_z,b_t)$. Any
restrictions on $F$ are explained in detail here.

 Then I specify the weights $(a_x,a_y,a_z,a_t)$ for a weighted blow up
and write down  the equation of the birational transform of $X$ on one
of the charts on the weighted blow up.  Before taking quotients, this
has the form
$t_1^{-m}F(x_1t_1^{a_x}, y_1t_1^{a_y}, z_1t_1^{a_z}, t_1^{a_t})$ if I
use the
$t$-chart. This is denoted by $B\tilde X$.

I need to take  quotient by 2 actions.  First   is the
$\frac1{a_t}(-a_x,-a_y,-a_z,1)$-action coming from the weighted blow
up. Second, the $\frac1{r}(b_x,b_y,b_z,b_t)$-action needs to be lifted
to the $(x_1,y_1,z_1,t_1)$-space. In some cases this lifts as a
$\z_r$-action but in other cases the actions combine into a
$\z_{(ra_t)}$-action. The quotient of $B\tilde X$ by these 2 actions
is a chart on the weighted blow up of $X$; it
 is denoted by $BX$. 

All these can be done in 4 different charts.  I chose the chart where
the singularities are most visible or the discrepancy computation is
the clearest.

Finally  I compute the exceptional divisor of the blow up, the
singularities of $BX$ and the discrepancy of the exceptional divisor.
\end{notation}

\begin{say}[$cAx/2$]

{\ } \newline  Normal form: $ax^2+by^2+g_{\geq 4}(z,t)/{\textstyle
\frac12(1,0,1,1)}$, where $ab\neq 0$.

Weights for blow-up: (1,1,1,1)

$t$-chart: $ax_1^2+by_1^2+t_1^{-2}g_{\geq 4}(z_1t_1,t_1)/{\textstyle
\frac12(0,1,0,1)}$

Exceptional divisor: $(t_1=ax_1^2+by_1^2=0)$. Over $\bar K$ this is
reducible and the two irreducible components are
$(t_1=\sqrt{a}x_1\pm \sqrt{-b}y_1=0)$. The $\z_2$-action interchanges
these two, so on the quotient we get a geometrically irreducible
exceptional divisor.

Singularity: The $\z_2$-action has a fixed curve on
$B\tilde X$: the intersection with the $(x_1=z_1=0)$-plane. Thus we
get a curve of nonterminal singularities on $BX$.

Discrepancy: $\pi^*\frac{dy\wedge dz\wedge dt}{x}= t_1\frac{dy_1\wedge
dz_1\wedge dt_1}{x_1}$, so $a(E,X)=1$.
\end{say}

\begin{say}[$cAx/4$]

{\ } \newline  Normal form: $ax^2+by^2+g_{\geq 2}(z,t)/{\textstyle
\frac14(1,3,1,2)}$, where $ab\neq 0$ and $g_2(0,1)=0$ for weight
reasons.

Weights for blow-up: (1,1,1,1)

$t$-chart: $ax_1^2+by_1^2+t_1^{-2}g_{\geq 2}(z_1t_1,t_1)/{\textstyle
\frac14(3,1,3,2)}$

Exceptional divisor: $\tilde E:=(t_1=ax_1^2+by_1^2+g_2(1,0)z_1^2=0)$. 
$\tilde E$ is geometrically irreducible if $g_2(1,0)\neq 0$. If
$g_2(1,0)= 0$, then $\tilde E$  is reducible over $\bar K$, and the
two irreducible components are
$(t_1=\sqrt{a}x_1\pm \sqrt{-b}y_1=0)$. The $\z_4$-action interchanges
these two, so on the quotient we get a geometrically irreducible
exceptional divisor $E$.

Singularity: The origin is on $B\tilde X$  since
$g_2(0,1)=0$ and it is a fixed point. We get an index 4 point on $BX$

Discrepancy: $\pi^*\frac{dy\wedge dz\wedge dt}{x}= t_1\frac{dy_1\wedge
dz_1\wedge dt_1}{x_1}$, so $a(E,X)=1$.
\end{say}

\begin{say}[$cD_4$ main series]

{\ } \newline  Normal form: $x^2+f_{\geq 3}(y,z,t)$, where we assume
that  $f_3(y,z,t)$ is irreducible over $\bar K$. 

Weights for blow-up: (2,1,1,1)

$x$-chart: $x_1+x_1^{-3}f_{\geq 3}(y_1x_1, z_1x_1,t_1x_1)/{\textstyle
\frac12(1,1,1,1)}$

Exceptional divisor: $\tilde E:=(x_1=f_3(y_1,z_1,t_1)=0)$. 
$\tilde E$ is geometrically irreducible by our assumption.

Singularity: The origin is a fixed point on $B\tilde X$, hence we get
an index 2 point on $BX$

Discrepancy: $\pi^*\frac{dy\wedge dz\wedge dt}{x}= x_1\cdot dy_1\wedge
dz_1\wedge dt_1$, so $a(E,X)=1$.
\end{say}

\begin{say}[$cD_4/2$ main series]

{\ } \newline  Normal form: $x^2+f_{\geq 3}(y,z,t)/{\textstyle
\frac12(1,1,0,1)}$, where we assume that  $f_3(y,z,t)$ is irreducible
over $\bar K$.  However, for weight reasons
$z|f_3(y,z,t)$, so this can not happen.
\end{say}

\begin{say}[$cD/3$]

{\ } \newline  Normal form: $x^2+f_{\geq 3}(y,z,t)/{\textstyle
\frac13(0,1,1,2)}$, where  $f_3(0,0,t)\neq 0$. Since this is not a
$cE$ point and for weight reasons, also $f_3(y,z,0)\neq 0$.  We can
write $f_3=t^3+f_3(y,z,0)$.

Weights for blow-up: (2,1,1,1)

$x$-chart: $x_1+x_1^{-3}f_{\geq 3}(y_1x_1, z_1x_1,t_1x_1)/{\textstyle
\frac12(1,1,1,1)}$, and then take the $\z_3$-action.

Lifting of the $\z_3$-action: It lifts to 
$\frac16(3,5,5,1)$. 
 
Exceptional divisor:  $\tilde E:=(x_1=f_3(y_1,z_1,t_1)=0)$. 
$\tilde E$ is geometrically irreducible  if $f_3(y,z,0)$ is not a
cube. If $f_3(y,z,0)=-L(y,z)^3$ over
$\bar K$, then $\tilde E$ has three geometrically irreducible
components
$(x_1=t_1-\eta L(y_1,z_1)=0)$ where $\eta^3=1$. The
$\z_6$-action permutes these, so on $BX$ we get a  geometrically
irreducible exceptional divisor.

Singularity: The origin is a $\z_6$-fixed point which has multiplicity
1 on $B\tilde X$.  $BX$ has a  terminal quotient singularity of index
6. 

Discrepancy: $\pi^*\frac{dy\wedge dz\wedge dt}{x}= x_1\cdot dy_1\wedge
dz_1\wedge dt_1$, so $a(E,X)=1$.

\end{say}

\begin{say}[$cD_{>4}$ and special $cD_4$]

{\ } \newline  Normal form: $x^2+Q_2(y,z,t)z+g_{\geq 4}(y,z,t)$, where
$Q_2(y,0,t)\neq 0$.  In the $cD_{>4}$ we always have this form (with
$Q_2(y,z,t)=y^2$). In the
$cD_4$-case we can achieve this form iff $f_3(y,z,t)$ has a simple
linear factor over $K$. 

Weights for blow-up: (2,1,2,1)

$z$-chart: $x_1^2+Q_2(y_1,z_1,t_1)+z_1^{-4}g_{\geq 4}(y_1z_1,
z_1^2,t_1z_1)/{\textstyle
\frac12(0,1,1,1)}$

Exceptional divisor: $\tilde
E:=(z_1=x_1^2+Q_2(y_1,0,t_1)+g_4(y_1,0,t_1)=0)$.  $\tilde E$ is
geometrically irreducible iff $Q_2(y_1,0,t_1)$ is not a square over
$\bar K$ or
$g_4(y_1,0,t_1)\neq 0$. If
$Q_2(y_1,0,t_1)=-L_1(y_1,t_1)^2$ (over $\bar K$) and
$g_4(y_1,0,t_1)=0$ then $\tilde E$ is   reducible over $\bar K$, and
the two irreducible components are
$(z_1=x_1\pm L_1(y_1,t_1)=0)$. The $\z_2$-action interchanges these
two, so $E\subset BX$ is  geometrically irreducible.

Singularity: The origin is a fixed point on $B\tilde X$, hence we get
an index 2 point on $BX$

Discrepancy: $\pi^*\frac{dy\wedge dz\wedge dt}{x}=
2z_1\frac{dy_1\wedge dz_1\wedge dt_1}{x_1}$, so $a(E,X)=1$.
\end{say}

\begin{say}[$cD_{>4}/2$ and special $cD_4/2$]

{\ } \newline  Normal form: $x^2+Q_2(y,z,t)z+g_{\geq
4}(y,z,t)/{\textstyle
\frac12(1,1,0,1)}$.

Weights for blow-up: (2,1,2,1)

$z$-chart: $x_1^2+Q_2(y_1,z_1,t_1)+z_1^{-4}g_{\geq 4}(y_1z_1,
z_1^2,t_1z_1)/{\textstyle
\frac12(0,1,1,1)}$ and then take the $\z_2$-action.

Lifting of the $\z_2$-action: We get a pair of commuting
$\z_2$-action on
$B\tilde X$, given by
$\frac12(0,1,1,1)$ and $\frac12(1,0,1,0)$. 

Singularity: The second action has a fixed curve on $B\tilde X$, so 
$BX$ is singular along a  curve.

Exceptional divisor and discrepancy: as in the $cD_{>4}$-case.
\end{say}

\begin{say}[$cD$-cases, conclusion]

{\ } \newline  We have settled all the $cD_{>4}$, $cD/2$ and
$cD/3$  cases, they have no g--extractions.

 In the $cD_4$ cases there are no g--extractions if $f_3$ is
irreducible or if it has a simple  linear factor over
$K$. The only remaining case is when $f_3$ is the product of 3 linear
factors which are conjugate over $K$.

This can not happen when $K=\r$, so over $\r$  points of type 
$cD, cD/2$ and $cD/3$  do not have g--extractions.

The situation is more complicated over  fields which do have cubic
extensions, as the following example shows. I have not classified all
cases.
\end{say}

\begin{exmp} Consider $x^2+y^3+az^3+t^6$, where $a\in K$ is not a
cube. The exceptional divisor of the $(3,2,2,1)$-blow up is
irreducible and has discrepancy 1. It has three points of index 2 
which are conjugate over $K$, and no other singularities. Hence this
is a g--extraction.
\end{exmp}

\begin{exmp} We obtain an interesting example from the equation
$x^2+(y^2+z^2)z+t^5$. The $(2,1,2,1)$ blow up has terminal
singularities (one with index 2). The exceptional divisor $E$  is
singular along a curve.
\end{exmp}

\section{$cA$-type Points}

In this section we study g--extractions of $cA$ type terminal
singularities. The conventions of (\ref{7.notation}),  (\ref{8.wbup})
and of (\ref{8.notation}) are used throughout.

\begin{say}[$cA_0$] (That is, smooth points.) {\ } \newline  Normal
form: ${\Bbb A}^3$.

The blow up of  the origin is  smooth with exceptional divisor
$E\cong \p^2$. $a(E,X)=2$, and by (\ref{gw.d-ineq.prop})
$a(F,X)\geq 2$ for every exceptional divisor
$F$ with $\cent_XF=\{0\}$. Therefore by (\ref{gw.discr2.cor}), the
blow up of  the origin is the only g--extraction. 

The exceptional divisor is $E\cong \p^2$ with normal bundle
$\o_{\p^2}(-1)$. 
\end{say}

\begin{say}[$cA_0/n,\ n\geq 2$]

{\ } \newline  Normal form: ${\Bbb A}^3/{\textstyle
\frac1{n}(r,-r,1)}$, where $(r,n)=1$ and $1\leq r\leq n-1$. 

Weights for blow-up: (r,n-r,1). 

$x$-chart: ${\Bbb A}^3(x_1,y_1,z_1)/{\textstyle
\frac1{r}(1,-n,-1)}$.

Exceptional divisor: $\tilde E:=(x_1=0)$.  Geometrically irreducible
and invariant under the
$\z_r$-action.

Lifting of the $\z_n$-action: The $\z_n$-action lifts to
$\frac1{n}(1,0,0)$. Its invariants are $x_2:=x_1^n$ and
$y_1,z_1$. The $\z_r$-action descends to the quotient of the
$\z_n$-action as
${\Bbb A}^3(x_2,y_1,z_1)/{\textstyle
\frac1{r}(n,-n,-1)}$. 

Singularity:  We obtain an index $r$ point on  the
$x$-chart, and similarly an index $n-r$ point on  the
$y$-chart.

Discrepancy: $\pi^*dx\wedge dy\wedge dz= rx_1^ndx_1\wedge dy_1\wedge
dz_1=
\frac{r}{n}x_1dx_2\wedge dy_1\wedge dz_1$. Since
$x_1=x_2^{1/n}$, we obtain that 
$a(E,X)=1/n$.

Conclusion: The above blow up is the only possible g--extraction. If
$n\geq 3$ then either $r\geq 2$ or
$n-r\geq 2$, and we obtain a singular point of index
$\geq 2$ on $BX$.

If $r=2$ then $BX$ is smooth, the exceptional divisor is
$E\cong \p^2$ with normal bundle $\o_{\p^2}(-2)$. $BX\to X$ is the
unique g--extraction.
\end{say}

\begin{say}[$cA_1$]

{\ } \newline  Normal form: $ax^2+by^2+cz^2 +dt^m$, where
$abcd\neq 0$.

Weights for blow-up: (1,1,1,1)

$t$-chart: $ax_1^2+by_1^2+cz_1^2 +dt_1^{m-2}$.

Exceptional divisor: $E:=(t_1=ax_1^2+by_1^2+cz_1^2=0)$ for
$m\geq 3$ and
$(t_1=ax_1^2+by_1^2+cz_1^2+d=0)$ for $m=2$. 
$E$ is geometrically irreducible.

Singularity: $BX$ has exactly one singular point for
$m\geq 4$, it lies on the
$t$-chart. $BX$ is smooth for $m=2,3$.

Discrepancy:  $\pi^*\frac{dy\wedge dz\wedge dt}{x}=
t_1\frac{dy_1\wedge dz_1\wedge dt_1}{x_1}$, so $a(E,X)=1$.

Conclusion: The only g--extraction is this blow up. The singularities
can be resolved by repeatedly blowing up the unique singular point.

\end{say}

\begin{say}[$cA_1/2$]

{\ } \newline  Normal form: $ax^2+by^2+cz^n +dt^m/{\textstyle
\frac12(1,1,1,0)}$, where
$abcd\neq 0$ and $\min\{n,m\}=2$.

Weights for blow-up: (1,1,1,1)

$z$-chart: $ax_1^2+by_1^2+cz_1^{n-2} +dt_1^mz_1^{m-2}$.

Exceptional divisor: $\tilde E:=(z_1=ax_1^2+by_1^2+c=0)$ for
$m\geq 3$,
$(z_1=ax_1^2+by_1^2+dt_1^2=0)$ for $n\geq 3$ and
$(z_1=ax_1^2+by_1^2+c+dt_1^2=0)$ for $n=m=2$. $E$ is geometrically
irreducible.

Singularity: The $\frac12(1,1,1,0)$ action lifts to a
$\frac12(0,0,1,1)$  action. Thus we get a fixed curve, where the blow
up intersects the plane $(z_1=t_1=0)$.

Discrepancy: $\pi^*\frac{dy\wedge dz\wedge dt}{x}= z_1\frac{dy_1\wedge
dz_1\wedge dt_1}{x_1}$, so $a(E,X)=1$.

Conclusion: The only possible g--extraction is this blow up. It has
nonterminal  singularities, so this does not occur.

\end{say}

\begin{say}[$cA_1/n,\ n\geq 3$]

{\ } \newline  Normal form: $xy+cz^{pm}+dt^2/{\textstyle
\frac1{n}(r,-r,1,0)}$, where $(r,n)=1$ and $cd\neq 0$. 

Weights for blow-up: (1,1,1,1)

$z$-chart: $x_1y_1+cz_1^{pm-2}+dt_1^2/{\textstyle
\frac1{n}(r-1,1-r,1,-1)}$

Exceptional divisor: $\tilde E:=(z_1=x_1y_1+dt_1^2=0)$, it is
geometrically irreducible.

Singularity: The $\z_n$-action has an isolated fixed  point at the
origin on $B\tilde X$. Thus $BX$ has an index
$n$ point.

Discrepancy: $\pi^*\frac{dy\wedge dz\wedge dt}{x}= z_1\frac{dy_1\wedge
dz_1\wedge dt_1}{x_1}$, so $a(E,X)=1$.

Conclusion: The only possible g--extraction is this blow up. It has a
higher index point, so this does not occur.
\end{say}

\begin{say}[$cA_{>1}^-$]

{\ } \newline  Normal form: $xy+g_{\geq 3}(z,t)$.

Weights for blow-up: (1,1,1,1)

$t$-chart: $x_1y_1+t_1^{-2}g_{\geq 3}(z_1t_1,t_1)$.

Exceptional divisor: $E:=(t_1=x_1y_1=0)$.  It has two  geometrically
irreducible  components.

Singularity: Not important

Discrepancy: $\pi^*\frac{dy\wedge dz\wedge dt}{x}= t_1\frac{dy_1\wedge
dz_1\wedge dt_1}{x_1}$, so $a(E,X)=1$.

Conclusion: At least 2 geometrically irreducible  divisors with
discrepancy
$\leq 1$. 

\end{say}

\begin{say}[$cA_{>1}/n,\ n\geq 3$ and $cA_{>1}^-/2$]

{\ } \newline  Normal form: $xy+g_{\geq 3}(z,t)/{\textstyle
\frac1{n}(r,-r,1,0)}$, where $(r,n)=1$. 

Weights for blow-up: (1,1,1,1)

$t$-chart: $x_1y_1+t_1^{-2}g_{\geq 3}(z_1t_1,t_1)/{\textstyle
\frac1{n}(r,-r,1,0)}$

Exceptional divisor: $\tilde E:=(t_1=x_1y_1=0)$. It is reducible and
both irreducible components are geometrically irreducible and
invariant under the
$\z_n$-action.

Singularity: Not important

Discrepancy: $\pi^*\frac{dy\wedge dz\wedge dt}{x}= z_1\frac{dy_1\wedge
dz_1\wedge dt_1}{x_1}$, so $a(E,X)=1$.

Conclusion: At least 2 geometrically irreducible  divisors with
discrepancy
$\leq 1$. 
\end{say}

\begin{say}[$cA_{>1}^+$, $\mult_0g$ even]

{\ } \newline  Normal form: $ax^2+by^2+g_{\geq 2m}(z,t)$, where $m\geq
2$, $-ab$ is not a square and
$g_{2m}\neq 0$. 

Weights for blow-up: (m,m,1,1)

$t$-chart: $ax_1^2+by_1^2+t_1^{-2m}g_{\geq 2m}(z_1t_1,t_1)$.

Exceptional divisor:
$E:=(t_1=ax_1^2+by_1^2+g_{2m}(z_1,1)=0)$, it is geometrically
irreducible.

Singularity: The $t$-chart on $BX$  is singular only at points $P$
corresponding to the multiple roots of
$g_{2m}(z,1)$. The singularity at $P$ again has type
$cA_{>1}^+$, but the multiplicity of the corresponding
$g^P(z_1,t_1)$ is not necessarily even. The
$z$ chart is similar.

Discrepancy: $\pi^*\frac{dy\wedge dz\wedge dt}{x}= t_1\frac{dy_1\wedge
dz_1\wedge dt_1}{x_1}$, so $a(E,X)=1$.

$x$-chart: $a+by_1^2+x_1^{-2m}g_{\geq 2m}(z_1x_1,t_1x_1)/ {\textstyle
\frac1{m}(1,0,-1,-1)}$.

Singularity: The fixed points of the $\z_m$-action are along the
$y_1$-axis, this itersects $B\tilde X$ in two points
$(0,\sqrt{-a/b},0,0)$ which are conjugate over $K$. Thus
$BX$ has 2  index $m$ terminal singularities  which are conjugate over
$K$.
 No other new singular points. The $y$-chart is similar.

Conclusion: The only g--extraction is the above weighted blow up.
 The exceptional divisor  is geometrically irreducible with a pair of
conjugate index
$m$-points. The other singular $K$-points of $BX$ are again  of type
$cA_{>1}^+$ or $cA_1$. 
\end{say}

\begin{say}[$cA^+_{>1}/2$, $\mult_0g$ even] 

{\ } \newline  Normal form: $ax^2+by^2+ g_{\geq 2m}(z,t)/{\textstyle
\frac12(1,1,1,0)}$, where $-ab$ is not a square and
$g_{2m}\neq 0$. 

Weights for blow up: $(m,m,1,1)$. 

$z$-chart: $ax_1^2+by_1^2+z_1^{-2m}g_{\geq 2m}(z_1,t_1z_1)
/{\textstyle \frac12(1-m,1-m,1,1)}$.

Exceptional divisor: $\tilde E:=(z_1=ax_1^2+by_1^2+g_{2m}(1,t_1)=0)$
is geometrically irreducible.

Singularities: If $m$ is odd then $(z_1=t_1=0)$  intersects $B\tilde
X$ in  a fixed curve of the
$\z_2$-action, thus we get a singular curve on $BX$.

If  $m$ is even, then on the $z$-chart the only
$\z_2$-fixed point is the origin. This is not on $B\tilde X$ iff
$z^{2m}\in g$. 

$t$-chart: $ax_1^2+by_1^2+t_1^{-2m}g_{\geq 2m}(z_1t_1, t_1)
/{\textstyle \frac12(1,1,1,0)}$.

 Singularities: On the $t$-chart the   fixed point set is the
$t_1$-axis. This intersects the exceptional divisor at the origin.
 This is not on the blow up iff $t^{2m}\in g$. 

Discrepancy: $\pi^*\frac{dy\wedge dz\wedge dt}{x}= t_1\frac{dy_1\wedge
dz_1\wedge dt_1}{x_1}$, so $a(E,X)=1$.

$x$-chart: $a+by_1^2+x_1^{-2m}g_{\geq 2m}(z_1x_1, t_1x_1) /{\textstyle
\frac1{m}(1,0,-1,-1)}$ and  we also need to take the quotient by the
$\z_2$-action.

Lifting the $\z_2$-action: $a+by_1^2+x_1^{-2m}g_{\geq 2m}(z_1x_1,
t_1x_1) /{\textstyle \frac1{2m}(1,0,m-1,-1)}$

 Singularities: On the $x$-chart the   fixed point set is the
$y_1$-axis. This intersects $B\tilde X$ at two points
$(0,\pm\sqrt{-a/b},0,0)$. We get a conjugate pair of terminal
singularities of index $2m$ on $BX$.

$y$-chart: Similar to  the $x$-chart.

Conclusion: $ax^2+by^2 +g_{\geq 2m}(z,t)/{\textstyle
\frac12(1,1,1,0)}$  where $-ab$ is not a square has a g--extraction
iff $m$ is even and $z^{2m},t^{2m}\in g_{2m}(z,t)$.  Under these
assumptions, the unique g--extraction is the $(m,m,1,1)$-blow up.
\end{say}

\begin{say}[$cA_{>1}^+$, $\mult_0g$ odd]\label{cA+.multg-odd}

{\ } \newline  Normal form: $ax^2+by^2+g_{\geq 2m+1}(z,t)$, where
$m\geq 1$, $-ab$ is not a square and
$g_{2m+1}\neq 0$. 

Weights for blow-up: $(s,s,1,1)$ for $1\leq s\leq m$,
 giving $B_sX\to X$.

$t$-chart: $ax_1^2+by_1^2+t_1^{-2s}g_{\geq 2m+1}(z_1t_1,t_1)$.

Exceptional divisor:
$E:=(t_1=ax_1^2+by_1^2=0)$, it is irreducible over $K$ but 
geometrically reducible.

Singularity:  The exceptional divisor itself has only smooth or normal
crossing  points, thus $B_sX$ has only
$cA$ type points. The
$(x_1=y_1=0)$ line is singular if $s<m$ and generically smooth for
$s=m$.
$B_mX$ is terminal. 

Discrepancy: $\pi^*\frac{dy\wedge dz\wedge dt}{x}= t_1\frac{dy_1\wedge
dz_1\wedge dt_1}{x_1}$, so $a(E,X)=1$.

Divisors with discrepancy 1:  Take the $(1,1,1,1)$-blow up. 
$BX$ is singular along a line with an $A_{2m-2}$ transversal section.
We can blow up the line
$(m-1)$-times. At each time the exceptional divisor is a pair of
transversally intersecting planes, thus we have only $cA$ type
singularities. After $(m-1)$ blow ups we obtain $g:Y\to X$ and $Y$ has
only isolated  $cA$ points, hence terminal. By
(\ref{gw.d-ineq.prop.cor}), all the exceptional divisors  over $0\in
X$ with discrepancy $1$ are birational to divisors on $Y$. They all
come in conjugate pairs and have been enumerated by the above
$(s,s,1,1)$ blow ups. 

Conclusion: There is  a unique g--extraction whose exceptional divisor
has discrepancy 1. It is the
$(m,m,1,1)$-blow up $B_mX\to X$. Its exceptional divisor is
geometrically reducible, so we need to look further.

\medskip

Divisors with discrepancy 2:  Let $F$ be a geometrically irreducible
exceptional divisor  over $0\in X$ with discrepancy $2$. Then
$\cent_{B_mX}F$ is real.  The center can not be the whole
$(x_1=y_1=0)$ line or a smooth point on it since both would give
$a(F,X)\geq 3$. Thus it is one of the singular points, corresponding
to a linear factor of $g_{2m+1}$.

By a linear change of the $z,t$-coordinates we may assume that this
linear factor is $z$. Thus $\cent_{B_mX}F$ is the origin of the $t$
chart, where $B_mX$ has equation 
$ax_1^2+by_1^2+t_1^{-2m}g_{\geq 2m+1}(z_1t_1,t_1)$. This is again a 
$cA_{>1}^+$ type point, ($\mult_0g$ can be even or odd) and
$a(F,B_mX)=1$.  We have already enumerated  all these cases, and we
know that $F$ is obtained by an
$(r,r,1,1)$-blow up. Putting the two steps together, we see that $F$
is obtained from $X$ by an
$(m+r,m+r,2,1)$-blow up. Next we   compute these.

 Normal form: $ax^2+by^2+g_{\geq 2m+1}(z,t)$, where
$m\geq 1$, $-ab$ is not a square,
$g_{2m+1}\neq 0$ and $\mult_0g(Z^2,T)\geq 2(m+r)$. 

Weights for blow-up: $(m+r,m+r,2,1)$, giving
 $B^rX\to X$.

$z$-chart: 
$ax_1^2+by_1^2+z_1^{-2(m+r)}g_{\geq 2m+1}(z_1^2,t_1z_1)/{\textstyle
\frac12(m+r,m+r,1,1)}$.

Singularity:  If $m+r$ is even, then the action has a fixed curve  on
$B^r\tilde X$, so $B^rX$ is not terminal. If $m+r$ is odd and 
 the origin is in $B^r\tilde X$, then  we get an index 2 point.
$z_1^{-2(m+r)}g_{\geq 2m+1}(z_1^2,t_1z_1)$ does not vanish at the
origin iff $z^{m+r}\in g_{\geq 2m+1}(z,t)$. This implies that $r\geq
m+1$.  But
$g_{2m+1}(z_1^2,t_1z_1)$ itself is not divisible by
$z_1^{4m+3}$, hence $r=m+1$.

Conclusion: Assume that there is a  linear change of the
$(z,t)$-coordinates  such that 
$$ g_{\geq 2m+1}(z,t)=\sum_{2i+j\geq 2m+2r}
\gamma_{ij}z^it^j,\qtq{and 
$\gamma_{ij}\neq 0$ for some $2i+j=2m+2r$.}
$$ In this coordinate system,  the  $(m+r,m+r,2,1)$ blow-up  gives  
an elementary extraction whose exceptional divisor is geometrically
irreducible and has discrepancy 2.
 Thus the only  possible g--extractions are  this weighted blow up and
the $(m,m,1,1)$ blow up found earlier.

The $(m+r,m+r,2,1)$ blow up  is  a g--extraction only in the $r=m+1$
case:
$$ g_{\geq 2m+1}(z,t)=\sum_{2i+j\geq 4m+2}
\gamma_{ij}z^it^j,\qtq{and} 
\gamma_{2m+1,0}\neq 0.
$$ In some cases  (cf. (\ref{cA.d=3.ex})), we do not have  any
geometrically irreducible exceptional divisor  over $0\in X$ with
discrepancy $2$. Then we have to compute further with discrepancy 3.
Fortunately, we can stop there.

\medskip

Divisors with discrepancy 3: 

 Normal form: $ax^2+by^2+g_{\geq 2m+1}(z,t)$, where
$m\geq 1$, $-ab$ is not a square and
$g_{2m+1}\neq 0$. 

Weights for blow-up: $(2m+1,2m+1,2,2)$, giving $Y\to X$.

$z$-chart: 
$ax_1^2+by_1^2+z_1^{-4m-2}g_{\geq 2m+1}(z_1^2,t_1z_1^2)/{\textstyle
\frac12(1,1,1,0)}$.

Exceptional divisor:
$E:=(t_1=ax_1^2+by_1+g_{2m+1}(1,t_1)=0)$, it is   geometrically
irreducible.

Singularity:  We get an index 2 point corresponding to the linear
factors of 
$g_{2m+1}(z,t)$. Thus over $\r$ there is always an index 2 point.

Discrepancy:   $\pi^*\frac{dy\wedge dz\wedge dt}{x}=
2z_1^3\frac{dy_1\wedge dz_1\wedge dt_1}{x_1}$, so $a(E,X)=3$.

\medskip

Final conclusion: These singularities always have g--extractions. One
is the
$(m,m,1,1)$-blow up. Its exceptional divisor is geometrically
reducible. This is the only g--extraction with discrepancy 1.

In some cases after a suitable coordinate change we can also perform
the
$(2m+1,2m+1,2,1)$ blow up. This is  the only g--extraction whose
exceptional divisor is  geometrically irreducible and has discrepancy
2.

If $g_{2m+1}(z,t)$ has no linear factors over $K$, then the
$(2m+1,2m+1,2,2 )$ blow up is a g--extraction  whose exceptional
divisor is  geometrically irreducible and has discrepancy 3. This is
the only one such. This  case never happens  over $\r$.

There may be other g--extractions  whose exceptional divisor is 
geometrically reducible and has discrepancy 2. I have no such examples.
\end{say}

\begin{exmp}\label{cA.d=3.ex}  Consider the singularity
$X:=(x^2+y^2+z^m+t^n=0)$ for $m,n$ odd and $m+2\leq n\leq 2m-1$. The
above computations show that  there is no geometrically irreducible
exceptional divisor  over $0\in X$ with discrepancy $\leq 2$.
\end{exmp}

\begin{say}[$cA^+_{>1}/2$, $\mult_0g$ odd]  {\ } \newline   Normal
form: $ax^2+by^2+ g_{\geq 2m+1}(z,t)/{\textstyle
\frac12(1,1,0,1)}$, where $-ab$ is not a square and
$g_{2m+1}\neq 0$. 

Weights for blow up: For weight reasons, only even powers of $t$
appear in
$g$. Thus we can define an integer $r$ by
$2m+2r=\mult_0g(Z^2,T)$.
$r\leq m+1$ since $g_{2m+1}\neq 0$. We consider the 
$(m+r,m+r,2,1)$ blow up.

$z$-chart: $ax_1^2+by_1^2+z_1^{-2m-2r}g_{\geq 2m+1}(z_1^2,t_1z_1)
/{\textstyle \frac12(m+r,m+r,1,1)}$ and then we have to take the
quotient by the
$\frac12(1,1,0,1)$-action. This lifts to a 
$\frac12(1,1,0,1)$-action on $B\tilde X$. We get a pair of commuting
$\z_2$-actions.

Exceptional divisor: $\tilde E:=(z_1=ax_1^2+by_1^2+\sum_{2i+j=2m+2r}
\gamma_{ij}t_1^j=0)$ is geometrically irreducible.

Discrepancy: $\pi^*\frac{dy\wedge dz\wedge dt}{x}=
2z_1^2\frac{dy_1\wedge dz_1\wedge dt_1}{x_1}$, so
$a(\tilde E,X)=2$.

Singularities:  If $m+r$ is even then the $\z_2\times
\z_2$-action is free in codimension one. One of the elements acts by
$(0,0,1,1)$, thus we get a singular curve in $BX$.

If $m+r$ is odd then one of the elements acts by
$(0,0,1,0)$. Coordinates on the quotient are given by 
$x_1,t_1, z_2=z_1^2,t_1$  and we get the equation
$$ ax_1^2+by_1^2+z_2^{-m-r}h_{\geq m+r}(z_2,t_1^2z_2) /{\textstyle
\frac12(1,1,0,1)}
$$ where $h(Z,T^2)=g(Z,T)$. At the origin we get a
$\z_2$-fixed point unless
$z^{m+r}\in g$. Thus $r\geq m+1$. On the other hand
$r\leq m+1$, thus $r=m+1$.  Computing the $t$-chart shows that  $BX$
has an index 2 point unless $t^{4m+2}\in g$. 

Discrepancy: From this we see that 
$a(E,X)=1/2$ if  $m+r$ is odd and $a(E,X)=2$ if  $m+r$ is even.

Conclusion: If $m+r$ is odd then a  g--extraction exists iff 
$$ g_{\geq 2m+1}(z,t)=\sum_{2i+j\geq 4m+2}
\gamma_{ij}z^it^j,\qtq{and} 
\gamma_{2m+1,0}\neq 0 \neq \gamma_{0,4m+2}.
$$ If this holds then  the  $(2m+1,2m+1,2,1)$ blow-up  is the unique
g--extractions. It has a   geometrically irreducible exceptional
divisor   with discrepancy $1/2$.

If $m+r$ is even then there may exist g--extractions with discrepancy
1/2 or 1. These can be determined by classifying all $\z_2$-invariant
divisors of disrepancy 1
 and   pointwise 
$\z_2$-fixed divisors of disrepancy 2 or 3 over $\tilde X$. The first 
task is easy, and we never get any g--extractions this way.  The
second task is harder and it seems to require separate consideration
of about a dozen cases; I have not done all of them.  Fortunately,
these singularities can be easily excluded in the main theorems using
topological considerations.
\end{say}

\section{$cE$-type Points}

In this section we study g--extractions of $cE$ type terminal
singularities. The conventions of (\ref{7.notation}),  (\ref{8.wbup})
and of (\ref{8.notation}) are used throughout.

\begin{say}[$cE_6$ main series] 

{\ } \newline  Normal form: 
$x^2+y^3+yg_{\geq 3}(z,t)+h_{\geq 4}(z,t)$. 

Weights for blow-up: (2,2,1,1)

$y$-chart:
$x_1^2+y_1^2+h_{4}(z_1,t_1)+y_1\Phi(y_1,z_1,t_1)/{\textstyle
\frac12(0,1,1,1)}$.

Exceptional divisor: 
$\tilde E:=(y_1=x_1^2+h_{4}(z_1,t_1)=0)$.
$\tilde E$  is geometrically irreducible iff  $h_4$ is not a square
over $\bar K$. If $-h_4$ is  a square over
$K$, then $\tilde E$ has 2 geometrically irreducible components. In
the other cases $\tilde E$ is irreducible over $K$ but reducible over
$\bar K$. Both of the components are fixed by the
$\z_2$-action, so the same 3 cases happen for $E$.

Singularity: The origin is a fixed point of the
$\z_2$-action which is on
$B\tilde X$. So  we get an index 2 point on $BX$.

Discrepancy: $\pi^*\frac{dy\wedge dz\wedge dt}{x}=
2y_1\frac{dy_1\wedge dz_1\wedge dt_1}{x_1}$, so
$a(E,X)=1$.

Conclusion:  If $-h_4$ is  a square over $K$ then there are 2
geometrically irreducible divisors with discrepancy 1, so no
g--extractions. If 
$h_4$ is not a square over $\bar K$ then we get an index 2 point, so
again there are no g--extractions.
\end{say}

\begin{say}[$cE/2$] 

{\ } \newline  Normal form: 
$x^2+y^3+yg_{\geq 3}(z,t)+h_{\geq 4}(z,t)/{\textstyle
\frac12(1,0,1,1)}$.  By weight considerations $g_3=0$ and $h_5=0$.  
$h_4\neq 0$ since otherwise we would not have a terminal point. This
is a $cE_6/2$ point.

Weights for blow-up: (2,2,1,1)

$y$-chart:
$x_1^2+y_1^2+h_{4}(z_1,t_1)+y_1\Phi(y_1,z_1,t_1)/{\textstyle
\frac12(0,1,1,1)}$.

Lifting of the $\z_2$-action.  The $\z_2$-action lifts to
$\frac12(1,1,0,0)$. Thus on $B\tilde X$ we have  two commuting
$\z_2$-actions.

Exceptional divisor: 
$\tilde E:=(y_1=x_1^2+h_{4}(z_1,t_1)=0)$. It is geometrically
irreducible iff  $h_4$ is not a square over $\bar K$. If
$h_4=-Q_2(z_1,t_1)^2$ then $\tilde E$ has 2 geometrically irreducible
components
$(y_1=x_1\pm Q_2(z_1,t_1)=0)$. The $\frac12(1,1,0,0)$ action
interchanges the 2 components, thus $E\subset BX$ is geometrically
irreducible.

Singularity: The $\frac12(1,1,0,0)$ action has a fixed curve, thus we
get a nonterminal singular curve on
$BX$.

Discrepancy: $\pi^*\frac{dy\wedge dz\wedge dt}{x}=
2y_1\frac{dy_1\wedge dz_1\wedge dt_1}{x_1}$, so
$a(E,X)=1$.
\end{say}

\begin{say}[$cE_7$ main series]\label{cE7.main.ser} 

{\ } \newline  Normal form: 
$x^2+y^3+yg_{\geq 3}(z,t)+h_{\geq 5}(z,t)$. 

Weights for blow-up: (3,2,1,1)

$x$-chart:
$x_1+y_1^3x_1+y_1g_{3}(z_1,t_1)+h_{5}(z_1,t_1)+x_1\Phi(y_1,z_1,t_1)/{\textstyle
\frac13(1,1,2,2)}$.

Exceptional divisor: 
$\tilde E:=(x_1=y_1g_{3}(z_1,t_1)+h_{5}(z_1,t_1)=0)$. It is
geometrically irreducible iff $g_3$ and $h_5$ have no common factors.

Singularity: The origin is a fixed point of the
$\z_3$-action which is on
$B\tilde X$. So  we get an index terminal 3 point on
$BX$. In fact, it is the index 3 terminal  point ${\Bbb
A}^3/\frac13(1,1,2)$.

Discrepancy: $\pi^*\frac{dy\wedge dz\wedge dt}{x}= 2x_1\cdot
dy_1\wedge dz_1\wedge dt_1$, so
$a(E,X)=1$.

Conclusion: If $g_3$ and $h_5$ have no common factors then $E$ is
irreducible and there are no g--extractions. 
\end{say}

\begin{say}[$cE$ with common linear factors]\label{cE.comm.l.f}

{\ } \newline  Normal form: 
$$ x^2+y^3+yzG_2(z,t)+z^2Q_2(z,t)+zH_4(z,t)+yg_{\geq 4}(z,t)+h_{\geq
6}(z,t).
$$ The following cases are of this form:

$cE_8$: $x^2+y^3+yg_{\geq 4}(z,t)+h_{\geq 5}(z,t)$, if
  $h_5$ has a   linear factor over $K$, which we can call $z$.

$cE_7$:
$x^2+y^3+yg_{\geq 3}(z,t)+h_{\geq 5}(z,t)$, if 
$g_3$ and $h_5$ have a common   linear factor over $K$, which we can
call $z$.

$cE_6$:
$x^2+y^3+yg_{\geq 3}(z,t)+h_{\geq 4}(z,t)$ if there is a   linear
factor over $K$, which we can call $z$, such that
$z^2|h_4,\ z|g_3$ and $z|h_5$.

Weights for blow-up: (3,2,2,1)

$z$-chart:
$$
\begin{array}{rl} x_1^2+y_1^3&+y_1G_2(0,t_1)+Q_2(0,t_1)+H_4(0,t_1)\\
&+y_1g_4(0,t_1)+h_6(0,t_1)+z_1\Phi(y_1,z_1,t_1)/{\textstyle
\frac12(1,0,1,1)}.
\end{array}
$$

Exceptional divisor $\tilde E$  is geometrically irreducible: 
$$ (z_1=
x_1^2+y_1^3+y_1G_2(0,t_1)+Q_2(0,t_1)+H_4(0,t_1)+y_1g_4(0,t_1)+h_6(0,t_1)=0).
$$

Singularity: The origin is a fixed point of the
$\z_2$-action which is on
$B\tilde X$. So  we get an index 2 point on $BX$.

Discrepancy: $\pi^*\frac{dy\wedge dz\wedge dt}{x}=
2z_1\frac{dy_1\wedge dz_1\wedge dt_1}{x_1}$, so
$a(E,X)=1$.
\end{say}

\begin{say}[$cE_6$ with $h_4$  a square] 

{\ } \newline  Normal form: 
$x^2+y^3+yg_{\geq 3}(z,t)+h_{\geq 4}(z,t)$. 

Weights for blow-up: (1,1,1,1)

$t$-chart:
$$ x_1^2+y_1^3t_1+y_1t_1^2g_3(z_1,1)+y_1t_1^3g_4(z_1,1)+
t_1^2h_{4}(z_1,1)+t_1^3h_5(z_1,1)+t_1^4\Phi(y_1,z_1,t_1).
$$ The $z$-chart is similar. 

Exceptional divisor: 
$E:=(t_1=x_1=0)$, and the scheme theoretic exceptional divisor is
$2E$.

Singularity:  On the $t$-chart the singular set  is the line
$L:=(x_1=y_1=t_1=0)$. We determine the singularities along this line.
For a fixed value
$z_1=b\in \bar K$ we get a $cA$-point if $h_4(b,1)\neq 0$. If
$h_4(z,1)$ has a simple root at $b$ then we get a
$cD$-point. If $h_4(z,1)$ has a multiple root at $b$ then we still get
a $cD$ point if $g_3(b,1)\neq 0$ and a
$cE$-point if $h_5(b,1)\neq 0$. 

Hence, if $b\in K$  and we do not have a cDV point, then $h_4$ has a
multiple linear factor which also divides $g_3$ and $h_5$. This case
was settled in (\ref{cE.comm.l.f}). Assuming that this is not the
case, we obtain that $BX$ has $cDV$ points along $L$.

The $z$-chart is similar and easy computations show that the $x$ and
$y$-charts are smooth along $E$.

Discrepancy: $\pi^*\frac{dy\wedge dz\wedge dt}{x}= t_1\frac{dy_1\wedge
dz_1\wedge dt_1}{x_1}$, so
$a(E,X)=2$.

First conclusion:  $BX$ is not a g--extraction since it has a singular
curve.
$E$ is geometrically irreducible and $a(E,X)=2$, thus if
$g:Z\to X$ is a g--extraction with exceptional divisor
$F$ then $a(F,X)=1$ by (\ref{gw.discr2.cor}).

Computations: Here we determine all divisors $F$ over
$0\in X$ with
$a(F,X)=1$.  If $\cent_{BX} F$ is not on $L$ then
$a(F,X)\geq 3$, and if $\cent_{BX} F$ is a point on $L$ then
$a(F,X)\geq 2$.  Thus if $a(F,X)=1$ then
$\cent_{BX}F=L$ and  $a(F,BX)=0$. 

Along $L$  the threefold $BX$ has  transversal type
$A_5$ whose singularity is resolved by blowing up the line 3-times. By
explicit computation we see that only the first of these produces an
exceptional divisor $F$ with
$a(F,X)=1$. This is the same divisor that we encountered in the
$(2,2,1,1)$-blow up and so it was already accounted for.

Final conclusion: There is no g--extraction except possibly when there
is a $b\in \bar K\setminus K$ such that
$(z-bt)^2|h_4, (z-bt)|g_3,(z-bt)|h_5$. In these cases the same
divisibilities hold if we replace $b$ by its conjugates over $K$. Thus
$b$ is  quadratic over $K$, a root of $Q_2(z,1)$. If $F$ is any
divisor over $0\in X$ with $a(F,X)=1$   then its center in $BX$ is 
$(z_1-b=0)\in L$ or  its conjugate. Thus $F$ is geometrically
reducible.
\end{say}

\begin{say}[$cE_6$ last case]\label{cE6.g--extr.exist}

{\ } \newline  Normal form: 
$$ x^2+y^3+cQ_2(z,t)^2+yL_1(z,t)Q_2(z,t)+C_3(z,t)Q_2(z,t)+ yg_{\geq
4}(z,t)+h_{\geq 6}(z,t),
$$ where $Q_2$ is a quadratic form which is irreducible over $K$  and
$-c$ is  not a square in $K$.  By a coordinate change as in
\cite[I.12.6]{AGV85} we can bring this to the simpler form
$$ x^2+y^3+cQ_2(z,t)^2+ yg_{\geq 4}(z,t)+h_{\geq 6}(z,t),
$$ though this is not important.

Normal form and topology over $\r$: We can choose $Q_2$ to be positive
definite and  diagonalize it. $-c\in \r$ is not a square, so we can
choose $c=1$. Thus we get the normal form
$$ x^2+y^3+(z^2+t^2)^2+ yg_{\geq 4}(z,t)+h_{\geq 6}(z,t).
$$ By (\cite[4.9]{rat1}) we obtain that $X(\r)$ is homeomorphic to
$\r^3$.

g--extractions: As we discussed above, all the  g--extractions of $X$
have geometrically reducible exceptional divisors.

Construction of g--extractions: It turns out that in these cases there
is a g--extraction. By  above remark  we do not need to know this for
certain to understand the topology  over $\r$, thus I only outline the
construction.

Basic constructions of toric geometry are used without reference; see
\cite{Fulton93}  for an introduction.

Over $\bar K$ we can bring the equation to the form
$$ x^2+y^3+z^2t^2+ yg_{\geq 4}(z,t)+h_{\geq 6}(z,t).
$$ Let $e_x,e_y,e_z,e_t$ be a basis of $\r^4$.  Consider the vectors
$w_z=\frac18(3,2,2,1)$ and $w_t=\frac18(3,2,1,2)$. These vectors give
a triangulation of the simplex with vertices $e_x,e_y,e_z,e_t$ where
the edges are
$$ (e_x,w_z), (e_x,w_t), (e_y,w_z), (e_y,w_t), (e_z,w_z), (e_t,w_t).
$$ Let us take the corresponding toric blow up. One can  check by a
rountine but tedious computation that all singularities of $BX$ are
terminal and we get two index 3 points on the chart corresponding to
the simplex
$(e_x,e_y,w_z,w_t)$. 

Note that the above construction is symmetric in $z$ and
$t$. Thus if we start with a quadratic form
$Q_2=z^2+qt^2$ and  introduce new coordinates
$z'=z+\sqrt{q}t$ and $t'=z-\sqrt{q}t$ then 
$Q_2=z't'$ and any blow up which is symmetric in $z',t'$ can be
transformed back to a  blow up of $X$ defined over $K$. We need to
check that the two index 3 points become conjugates over $K$, but this
is easy to see from the explicit equations.
\end{say}

\begin{say}[$cE_7$ with common nonlinear factor]

{\ } \newline  Normal form: 
$x^2+y^3+yg_{\geq 3}(z,t)+h_{\geq 5}(z,t)$, where we assume that the
greatest common divisor  of 
$g_3$ and $h_5$ is $K$-irreducible (and nonconstant). We write
$g_3=Q(z ,t)G(z,t)$ and
$h_5=Q(z,t)H(z,t)$. ($Q$ is allowed to be linear, though this case is
treated already.)

Weights for blow-up: (3,2,1,1)

$x$-chart:
$x_1+y_1^3x_1+y_1g_{3}(z_1,t_1)+h_{5}(z_1,t_1)+x_1\Phi(y_1,z_1,t_1)/{\textstyle
\frac13(1,1,2,2)}$.

Exceptional divisor: It has two irreducible components: 
\begin{eqnarray*}
\tilde E&:=&(x_1=y_1G(z_1,t_1)+H(z_1,t_1)=0),\qtq{and}\\
\tilde F&:=&(x_1=Q(z_1,t_1)=0). 
\end{eqnarray*}
$\tilde E$ 
 is geometrically irreducible, $\tilde F$ is irreducible but 
geometrically reducible if $Q$ is not linear.

Discrepancy: $\pi^*\frac{dy\wedge dz\wedge dt}{x}=
3x_1\frac{dy_1\wedge dz_1\wedge dt_1}{x_1}$, so
$a(E,X)=1=a(F,X)$. (The latter equality uses that $Q$ is not a
multiple factor.)

Further aim: We would like to construct a birational morphism $g:Z\to
X$ whose exceptional divisor corresponds to $E$, and  determine the
singularities of
$Z$. Thus in $BX$ we have to contract $F$. $F$ is not
$\q$-Cartier in $BX$ and $F$ can not be contracted in
$BX$. First we have to correct this problem.

Singularities of $BX$: I claim that $BX$ has only canonical
singularities. This   can be done 2 ways. One can compute each chart
explicitly, which is  rather  tedious.    I found it easier to use  a
degeneration argument as follows.  Let $F$ be the normal form of the
equation as above.  We may assume that
$g_3(1,0)=1$. Consider the substitution
$$ F(x,y,z,t)\mapsto 
\epsilon^{-24}F(\epsilon^{12}x,\epsilon^{8}y,\epsilon^{6}z,
\epsilon^{7}t).
$$  The exponents are chosen so that for $\epsilon\to 0$ the limit is
$X_0:=(x^2+y^3+yz^3=0)$.  The $(3,2,1,1)$-blow up $BX_0$ is easy to
compute. We find an index 3 terminal  point
${\Bbb A}^3/\frac13(1,1,2)$, a curve of  $cA$-points and  a curve of
$cE_7$-points corresponding to the
$t$-axis. Thus $BX$, as a small deformation of $BX_0$, has an index 3
point at the origin and some $cDV$ singularities. (These turn out to
be isolated points but we do not need this.) As in
(\ref{cE7.main.ser}) we see that the index 3 point is at the origin of
the $x$-chart and it is 
${\Bbb A}^3/\frac13(1,1,2)$. In particular it is
$\q$-factorial.

Small blow up: Let $p:Y\to BX$ be the blow up of $F$ in
$BX$. Let $F'\subset Y$ denote the birational transform of $F$. Away
from the index 3 point
$BX$ is locally isomorphic to $B\tilde X$.  
$\tilde F$ is defined by 2 equations $(x_1=Q(z_1,t_1)=0)$, thus
$p:Y\to BX$ is small and is an isomorphism at all points where $F$ is
$\q$-Cartier. The index 3 point is
$\q$-factorial, so $F$ is $\q$-Cartier there. Thus
$p:F'\to F$ is an isomorphism.

Contracting $F'$: $F$ is a cone over a $K$-irreducible  curve, hence
its cone of curves over $K$ is 1-dimensional. If $C\subset F'$ is a
general curve then
$(C\cdot K_Y)=(p(C)\cdot K_{BX})<0$ and
$(C\cdot F')=(p(C)\cdot F)<0$. Thus the curves in $F'$ generate a
$K_Y$-negative extremal ray of  $Y/X$, which can be contracted. We
obtain
$f: Y\to Z$ and $g:Z\to X$. 
$P:=f(F')$ is a $K$-point since $F'$ is connected.  

Conclusion: $g:Z\to X$ has a geometrically irreducible exceptional
divisor corresponding to $E$ and it has discrepancy 1. Furthermore, by
(\ref{gw.discr2.cor})   the index of
$P$ can not be one since $a(F,X)=1$. Hence  there are no
g--extractions. 
\end{say}

\begin{say}[Conclusion] The $cE_8$ case is settled if
$h_5(z,t)$ has a linear factor over $K$. This always holds if  $K=\r$,
hence at least in this case  there are no g--extractions. I do not
know what happens if $K\neq
\r$.

The $cE_7$ case is settled if $g_3(z,t)$ and 
$h_5(z,t)$ have no common factor, or if they have a common linear
factor over
$K$ or if they have a unique common factor over $K$. This accounts for
all the possibilities, hence there are no g--extractions.

The $cE_6$ case is settled if $h_4(z,t)$ is not a square over $\bar K$,
 if $-h_4(z,t)$ is a square over $K$ or if $h_4(z,t)$ is divisible by
the square of a linear form over $K$. In these cases there are no
g--extractions.

The remaining case is treated in (\ref{cE6.g--extr.exist})  and the
unique g--extraction is written down explicitly. For the applications
in this paper the existence is not crucial.
\end{say}

\begin{exmp}  Let $X$ be the $cE_7$ type singularity
$x^2+y^3+yg_3(z,t)+h_5(z,t)$, where $g_3$ and $h_5$ do not have a
common factor. It is not hard to see that
$X$ is an isolated singular point and its
$(3,2,1,1)$-blow up has only terminal singularities. As in
(\ref{cE7.main.ser}), the $y$ chart on the blow up gives the
exceptional divisor
$$ E=(g_3(z,t)+h_5(z,t)=0)/{\textstyle
\frac12(1,1,1,1)}.
$$ This gives  examples of extremal contractions whose exceptional
divisor $E$ has a quite complicated singularity along the
$(z=t=0)$-line.
\begin{enumerate}
\item $x^2+y^3+yz^3+t^5$. $E$ is singular along 
$(z=t=0)$, with a transversal singularity type
$z^3+t^5$, that is $E_8$.
\item $x^2+y^3+y(z-at)(z-bt)(z-ct)+t^5$. $E$ has  triple
selfintersection  along 
$z=t=0$.
\end{enumerate}
\end{exmp}

\section{Hyperbolic 3--manifolds}

The aim of this section is to show that every  hyperbolic 3--manifold
satisfies the conditions (\ref{int.no.cond}).

\begin{thm}\label{hyp.doesnotcont.thm} Let $M$ be a compact 
hyperbolic 3--manifold. Then $M$ does not contain any PL submanifold
 of   the following types:
\begin{enumerate}
\item $\r\p^2$
\item 1--sided $S^1\times S^1$
\item 1--sided Klein bottle.
\end{enumerate}
\end{thm}

 We use two facts about hyperbolic 3--manifolds. First, that their
universal cover is homeomorphic to $\r^3$. Second, that their
fundamental group does not contain a subgroup isomorphic to $\z^2$
(see, for instance, \cite[4.6]{Scott83}). 
\medskip

More generally, we see how these conditions fit in the framework of
Thurston's geometrization conjecture. This version was pointed out to
me by  Kapovich.

\begin{thm}\label{gen.doesnotcont.thm} Let $M$ be a compact 
  3--manifold.  Assume that $M=M_1\ \#\ \cdots \ \#\ M_k$, where
\begin{enumerate}
\item[(i)] each $M_i$ is aspherical, and 
\item[(ii)] the Seifert fibered part of the Jaco--Shalen--Johannson
decomposition of  $M_i$  is orientable.
\end{enumerate} Then $M$ does not contain any PL submanifold
 of   the following types:
\begin{enumerate}
\item $\r\p^2$
\item 1--sided $S^1\times S^1$
\item 1--sided Klein bottle.
\end{enumerate}
\end{thm}

We consider the 3 types of submanifolds separately. Condition
(\ref{int.no.cond}.1) is closely related to the notion of
$\p^2$-irreducibility (cf.\ \cite[p.88]{Hempel76}).

\begin{lem} Let $M$ be a  3--manifold  with universal cover $\tilde M$.
\begin{enumerate}
\item  If $M\sim M_1\ \#\ M_2$, then $M$  contains a 2--sided $\r\p^2$
iff one of the summands does.

\item  Assume that  $\tilde M$ is homeomorphic to
$\r^3$. Then $M$ does not contain an $\r\p^2$ and $M$ can not be
written as  a nontrivial connected sum.
\end{enumerate}
\end{lem}

Proof.  Assume that  $F\subset M$ is a 2--sided $\r\p^2$. We may
assume that $F$ is transversal to the gluing $S^2$. Thus $C=F\cap S^2$
is an embedded curve in $F$. Assume first that
$F$ has a connected component $C_1\subset C$  which is not null
homotopic in
$F$. Then $F$ is not orientable along $C_1$, and the same holds for
$M$ along $C_1$. But $M$ is orientable along $S^2$, a contradiction.

Take  any connected component $C_i\subset C$ such that
$C_i\subset S^2$ bounds a disc $D_i$ which is disjoint from $C$.
$C_i$ also bounds a disc $D'_i$ in $F$ (since it is null homotopic in
$F$). Thus we can change the embedding $\r\p^2\to M$ by replacing
$D'_i$ with
$D_i$ and then pushing it to one side.  The new embedding is still
2--sided.  Repeating if necessary, we eventally get an embedding which
is disjoint from $S^2$, proving (1). 

 $\r\p^2$ can not be embedded into $\r^3$ (cf. \cite[27.11]{GrHa81}),
thus the preimage of
$\r\p^2$ in $\r^3$ is a union of copies of $S^2$. Fix one of these and
call it $N$. By the Schoenflies theorem (cf.\ \cite[Sec. 17]{Moise77})
  $N$ bounds a 3--ball $B^3$.  At least one element of $\pi_1(M)$ maps
$N$ to itself. It can not map the inside of
$N$ to its outside since these are not homeomorphic. If it maps 
  $B^3$ to itself, then by the Borsuk--Ulam theorem (cf.\
\cite[23.20]{Fulton95}) we have a  covering transformation with a
fixed point, a contradiction. 

Assume that we have $S^2\sim N'\subset M$ and let $S^2\sim N\subset
\tilde M$ be one of the preimages. Then $N$ bounds a 3--ball  and so
does $N'$. 
\qed
\medskip

In order to study the conditions (\ref{int.no.cond}.2--3) we have to
distinguish two cases.

\begin{say}[Incompressible case]\label{12.incompr.say} 

Let $M$ be a compact 3--manifold and $S\subset M$ a compact 1--sided
torus or Klein bottle. Assume that $\pi_1(S)\into \pi_1(M)$. Let
$\partial U$ be  the boundary of a regular neighborhood of $S$. Then 
$\partial U$ is a 2--sided  torus or Klein bottle and 
$\pi_1(\partial U)\into \pi_1(S)\into \pi_1(M)$ is an injection. This
implies that $\partial U$ is incompressible in $M$ (cf.\
\cite[pp.88-89]{Hempel76}). Thus $U$ is one of the pieces of the
Jaco--Shalen--Johannson decomposition of $M$ (cf.\
\cite[p.483]{Scott83}). We have to be a little more careful since
$U$ is Seifert fibered, thus it may sit inside one of the Seifert
fibered components. 

 The fundamental group  of a hyperbolic 3--manifold does not contain a
subgroup isomorphic to $\z^2$ (see, for instance,
\cite[4.6]{Scott83}), hence the incompressible case does not happen
for hyperbolic 3--manifolds.
\end{say}

\begin{say}[Compressible case]\label{12.compr.say} 

In this case we show that $M$ can be written as a connected sum   with
a very special summand.

\begin{prop}\label{1-s.torus}
 Let $M$ be a compact 3--manifold. Then $M$ contains a 1--sided torus 
$T$ such that $\pi_1(T)\to \pi_1(M)$ is not an injection iff 
$M\sim N\ \#\ (S^1\tilde{\times}S^2)$ or $M\sim N\ \#\ (S^1
\times\r\p^2)$
\end{prop}

Proof. Let  $T\subset U\subset M$ be a regular neighborhood. Set
$V=M\setminus U$. Then $\partial U=\partial V\sim S^1\times S^1$. We
know that
$\pi_1(\partial U)$ injects into $\pi_1(U)$. If $\pi_1(\partial
U)\into \pi_1(V)$, then
$\pi_1(\partial U)\into \pi_1(U)\into \pi_1(M)$ 
 by Schreier's theorem (cf. \cite[IV.2.6]{Lyndon-Schupp77}).
$\pi_1(\partial U)$ is an index 2 subgroup of $\pi_1(T)$ an $\pi_1(T)$
is torsion free. Thus 
$\pi_1(T)\to \pi_1(M)$  is also an injection, a contradiction.

Therefore, by the Loop theorem (cf.\ \cite[4.2]{Hempel76}), there is an
embedding of the  disc $j:(B,\partial B)\into (V,\partial V)$ such that
the image of $j(\partial B)$ is not contractible in $\partial V$. 

Let us cut $V$ along $j(B)$ to get $W$. 
 The   boundary of $W$ is $\partial V$ cut along
$j(\partial B)$ (which is a cylinder) with two copies of  $B$ pasted to
the ends. That is, $\partial W\sim S^2$. Therefore $M$ is obtained by
pasting   $W$ to a 3--manifold (with boundary)
$K$, which is obtained from 
$U$ by attaching a 2--handle.

There are two cases corresponding to whether $j(\partial B)$ is a
primitive element of $\pi_1(U)\cong \z^2$  (hence $\pi_1(K)\cong \z$)
or is contained in $2\z^2$ (hence $\pi_1(K)\cong \z+\z_2$).\qed

\begin{prop}\label{1-s.kb}
 Let $M$ be a compact 3--manifold which does not contain a 2--sided
$\r\p^2$. Then $M$ contains a 1--sided Klein bottle
$K$ such that $\pi_1(K)\to \pi_1(M)$ is not an injection iff 
$M\sim N\ \#\ (S^1\tilde{\times}S^2)$ or $ M\sim N\ \#\ (\r\p^3\ \#\
\r\p^3)$.
\end{prop}

Proof. Let  $K\subset U\subset M$ be a regular neighborhood. As in the
proof of (\ref{1-s.torus}) we obtain
 an embedding of the  disc $j:(B,\partial B)\into (V,\partial V)$ such
that the image of $j(\partial B)$ is not contractible in $\partial V$. 
We again cut $V$ along $j(B)$ to get $W$.  Let $\partial V^*$ denote
$\partial V$ cut along
$j(\partial B)$. 

There are 3 cases to consider corresponding to what $\partial V^*$ is:
\begin{enumerate}
\item  ($\partial V^*$ is a cylinder). Then we obtain a connected sum
decomposition as in (\ref{1-s.torus}).
\item ($\partial V^*$ consists of two Moebius bands). Then $\partial W$
is two disjoint projective planes, hence $M$   contains a 2--sided
projective plane. This can not happen by assumption.
\item  ($\partial V^*$ is a Moebius band). In this case
$j(\partial B)$ is 1--sided in $\partial V$, thus $M$ is not orientable
along $j(\partial B)$. Then $j(\partial B)$ can not be the boundary of
an embedded disc.\qed
\end{enumerate}
\end{say}

\begin{rem} So far we have  excluded Seifert fiber spaces from
considerations.  Many Seifert fiber spaces do contain  1--sided tori
or Klein bottles.

If $p:M\to F$ is  a Seifert fiber space and $C\subset F$ a 1--sided
curve not passing through any critical value, then $p^{-1}(C)\subset
M$ is a 1--sided torus or Klein bottle. Another example can be
obtained as follows. Let  $x,x'\in F$ be two points such that the
fibers over   them have multiplicity 2. Let $I\subset F$ be a simple
path connecting $x$ and
$x'$. Then $p^{-1}(I)$ is   a 1--sided   Klein bottle.

It is not hard to see that if $T\subset M$ is a 
 1--sided torus or Klein bottle such that $p(T)$ is 1--dimensional
(these are called vertical) then $T$ is obtained  by one of the  
above constructions.

Assume now in addition that $M$ has a geometry modelled on ${\mathbb
H}^2\times \r$ (cf.\ \cite[p.459]{Scott83}). Then by
\cite[5.6]{Johannson79}, every 1--sided torus or Klein bottle in $M$ is
isotopic to a vertical one.

This way we obtain  many examples of nonorientable  Seifert fiber
spaces which satisfy the conditions (\ref{int.no.cond}).
\end{rem}

\noindent University of Utah, Salt Lake City UT 84112 

\begin{verbatim}kollar@math.utah.edu\end{verbatim}

\end{document}